%% file: publication.tex
\documentclass{article} 
\usepackage{amsfonts}
\usepackage{amsopn}
\usepackage{epsfig}
\usepackage{amsmath}
\begin{document}
\title{Degenerations of Surface Scrolls and the Gromov-Witten
  Invariants of Grassmannians}
\author{Izzet Coskun}
\maketitle
\noindent {\bf Abstract:} We describe an algorithm for computing certain
characteristic numbers of rational normal surface scrolls using
degenerations. As a corollary we obtain an efficient method for
computing the corresponding Gromov-Witten invariants of the
Grassmannians of lines.
\tableofcontents
\bigskip
\bigskip
\bigskip
\bigskip
\bigskip
\bigskip
\bigskip
\bigskip
\bigskip
\bigskip
\bigskip
\bigskip

\thispagestyle{empty} 


\pagestyle{plain}
\newpage
\setcounter{page}{1}

\newtheorem{theorem}{Theorem}[section]
\newtheorem{proposition}[theorem]{Proposition}
\newtheorem{lemma}[theorem]{Lemma}
\newtheorem{definition}[theorem]{Definition}
\newtheorem{remark}{Remark}[section]
\newtheorem{corollary}[theorem]{Corollary}

\renewcommand{\P}{\mathbb{P}}
\newcommand{\C}{\mathbb{C}}
\renewcommand{\O}{\mathcal{O}}
\newcommand{\G}{\mathbb{G}}
\renewcommand{\a}{\sum_{i=1}^a \alpha_i}
\renewcommand{\b}{\sum_{j=1}^b \beta_j} 
\newcommand{\X}{\mathcal{X}}
\newcommand{\Y}{\mathcal{Y}}
\newcommand{\Z}{\mathbb{Z}}
\newcommand{\M}{\overline{M}_{0,n}}
\newcommand{\F}{\mathbb{F}}

\section{Introduction}

This paper investigates the enumerative geometry of rational normal
surface scrolls in $\P^N$ using degenerations. We obtain an effective
algorithm for computing certain characteristic numbers of balanced
scrolls. Surface scrolls can be interpreted as curves in the
Grassmannian of lines $\G(1,N)$. Using our algorithm we calculate the
corresponding Gromov-Witten invariants of $\G(1,N)$.  We work over the
field of complex numbers $\C$. \smallskip

\noindent {\bf Motivation.} By the {\it characteristic number
  problem} we mean the problem of computing the number of varieties in
$\P^n$ of a given `type' [e.g. curves of degree $d$ and genus $g$]
that meet the `appropriate' number of general linear spaces so that
the expected dimension is zero. This problem has attracted a lot of
interest since the 19th century (see \cite{schubert:enumerative},
\cite{kleiman:schubert}). In the last decade, motivated by the work of
string theorists and Kontsevich, there has been significant progress
on the problem for curves (see \cite{Capharris:severi},
\cite{vakil:quartic}, \cite{vakil:rationalelliptic} for references).

In comparison the charcateristic numbers of higher dimensional
varieties are harder to compute, hence have received less attention.
In this paper we start a more systematic study of the characteristic
numbers of higher dimensional varieties using degenerations.  Here we
restrict our attention to rational surface scrolls, although most of
the techniques apply with little change to higher dimensional scrolls
\cite{coskun3:degenerations} and can be modified to apply to Del Pezzo
surfaces \cite{coskun2:degenerations}.

The enumerative geometry of scrolls is also attractive for its
connection to the Gromov-Witten theory of $\G(1,N)$. Localization
techniques and the associativity relations in the quantum cohomology
ring lead to recursive algorithms that compute the invariants, but
these algorithms are usually inefficient. For example, using FARSTA
\cite{farsta}, a computer program that computes Gromov-Witten
invariants from associativity relations, it takes several months to
determine cubic and quartic invariants of $\G(1,5)$. The algorithm we
prove here allows us to compute many of these invariants by hand (\S
8, 9).  \smallskip

\noindent {\bf Notation.} Let $\overline{M}_{0,n}(\G(1,N),d)$ denote
the Kontsevich moduli space of $n$-pointed genus 0 stable maps to
$\G(1,N)$ of Pl\"ucker degree $d$. Let $\mbox{Hilb}(\P^N,S^d)$ denote
the component of the Hilbert scheme whose general point corresponds to
a smooth rational normal surface scroll $S$ of degree d in $\P^N$.
\smallskip

\noindent {\bf Results.} The main results of this paper are the
following: \smallskip

$\bullet$ To calculate the characteristic numbers of scrolls, we
specialize the linear spaces meeting the scrolls to a general
hyperplane $H$. We prove that a general, non-degenerate, reducible
limit of balanced scrolls incident to the linear spaces consists of
the union of two balanced scrolls meeting along a line---provided that
the limit of the hyperplane sections in $H$ remains non-degenerate.
(\S6) The precise statements are given in Theorem \ref{main1} and
\ref{main2}. \smallskip

$\bullet$ By successively breaking the scrolls to smaller degree
scrolls, we obtain a recursive algorithm for computing the
characteristic numbers of balanced scrolls in $\P^N$ incident to
linear spaces of small dimension (\S 8). Theorem \ref{enumerative}
summarizes the result. \smallskip

\noindent {\bf Example:} For instance, the algorithm easily shows that
the number of scrolls of degree $n$ in $\P^{n+1}$ containing $n+5$
general points and meeting a general $n-3$ plane is $(n-1)(n-2)$
(\S5). \smallskip
  
$\bullet$ As a corollary, we obtain an efficient method for computing
the corresponding Gromov-Witten invariants of $\G(1,N)$.  The proof
also yields a method for computing some Gromov-Witten invariants of
$\F(0,1;N)$, the partial flag variety of pointed lines in $\P^N$ (\S
9). \smallskip

\noindent {\bf The method.} Our method is a degeneration method
inspired by \cite{Capharris:severi} and especially
\cite{vakil:rationalelliptic}.  The prototypical example answers the
question how many lines meet 4 general lines $l_1, \cdots l_4$ in
$\P^3$. If $l_1$ and $l_2$ lie in a plane $P$, then the answer is easy
to see. Let $q = l_1 \cap l_2$. The two solutions are the line in $P$
through $l_3 \cap P$ and $l_4 \cap P$ and the intersection of the two
planes $\overline{q l_3}$ and $\overline{q l_4}$.
\begin{figure}[h!]
\begin{center}
\input{lines.pstex_t}
\end{center}
\caption{Prototypical example of the degeneration method.}
\label{Figure 1}
\end{figure}
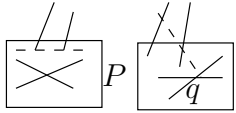
To answer the original question we can specialize two of the lines to
a plane. If we know how many of the original solutions approach each
of the two special solutions we can answer the problem. Our algorithm
carries out this classical idea for rational normal scrolls.

The solution of an enumerative problem by degenerations has two steps.
We specialize linear spaces meeting the scrolls to a general
hyperplane one at a time. First, we identify the limiting positions of
the scrolls. 

We prove that non-degenerate limits of scrolls are unions of scrolls
where any two adjacent components share a common fiber. These limit
surfaces occur as images of trees of Hirzebruch surfaces (Proposition
\ref{centralfiber}). We describe the trees that occur as limits of
scrolls $S_{k,l}$.

Not every tree of scrolls smooths to $S_{k,l}$. The specializations of
$S_{k,l}$ contain a connected degree $k$ curve whose components are
rational curves in section classes on the scrolls. The existence of a
degree $k$ curve with these properties turns out to be sufficient for
a union of two scrolls of total degree $k+l$ to smooth to $S_{k,l}$
(Proposition \ref{direc}).

In \S 6 we carry out a detailed dimension count to identify which
unions of scrolls occur as limits under some non-degeneracy
assumptions. The dimension calculations are considerably harder for
surfaces than for curves because the need to trace both the hyperplane
section and the directrix of the surface forces us to work with a
non-convex space.  Nonetheless, we prove that in a general
one-parameter family if the surfaces and their hyperplane sections
remain non-degenerate, balanced scrolls break into unions of balanced
scrolls.

Once we determine the limits, we need to determine their
multiplicities.  We reduce the calculations to the case of curves by
constructing a smooth morphism from the space of scrolls to
$\overline{M}_{0,n}(\P^N, d)$ and pulling-back the relations between
cycles in these spaces \cite{vakil:rationalelliptic} to the space of
scrolls (\S 7). We give many examples to illustrate how the algorithm
works in \S 5.  \smallskip

\noindent {\bf Acknowledgements.} I would like to thank J. de
Jong, M. Popa, J. Starr, D. Avritzer and M. Mirzakhani for fruitful
discussions. I am grateful to R. Vakil and my advisor J.  Harris for
invaluable suggestions on both the content and presentation of this
paper. I would also like to thank J. Li, the Stanford Mathematics
Department and especially R. Vakil for their hospitality.

\section{Preliminaries on Scrolls}

This section provides a summary of basic facts about rational scrolls
and systems of divisors on them; for more details consult
\cite{beauville:surface} Ch. 4, \cite{friedman:surfaces} Ch. 5
or \cite{joe:thebook} \S 3  Ch. 4.  \bigskip

\noindent {\bf Rational normal scrolls.} Let $k \leq l$ be two
non-negative integers with $l>0$. We will denote a rational normal
surface scroll of bidegree $k,l$ by $S_{k,l}$.  $S_{k,l}$ is a
rational surface of degree $k+l$ in $\P^{k+l+1}$. We now recall
its construction.  \smallskip

Fix two rational normal curves of degrees $k$ and $l$ in $\P^{k+l+1}$
with disjoint linear spans. Fix an isomorphism between the curves.
$S_{k,l}$ is the surface swept by the lines joining the points
corresponding under the isomorphism. The degree $k$ curve is called
the {\it directrix}. If the directrix reduces to a point, we obtain
$S_{0,l}$, the cone over a rational normal $l$ curve.  We will call a
scroll {\it balanced} if $l-k \leq 1$, and {\it perfectly balanced} if
$k=l$. A perfectly balanced scroll has a one-parameter family of
directrices. \smallskip

Rational normal scrolls are non-degenerate surfaces of minimal degree
in projective space.  Conversely, 
  
\begin{proposition}\label{classification}
  (\cite{joe:thebook} p.525) Every non-degenerate irreducible surface of
  degree $m-1$ in $\P^m$ is a rational normal scroll or the
  Veronese surface in $\P^5$.
\end{proposition}

\noindent {\bf Hirzebruch surfaces.} The {\it Hirzebruch surface}
$F_r$, $r \geq 0$, is the projectivization of the vector bundle
$\O_{\P^1} \oplus \O_{\P^1}(r)$ over $\P^1$. In this paper the
projectivization of a vector bundle $\P V$ will mean the one-dimensional
subspaces of $V$.  \smallskip

The Picard group of $F_r$, $r>0$, is generated by two classes: the
class $f$ of a fiber $F$ of the projective bundle and the
class $e$ of the unique section $E$ with negative
self-intersection. The intersection pairing is given by
\begin{displaymath}
f^2 = 0, \ \ \ \ f \cdot e = 1, \ \ \ \ e^2 = -r.
\end{displaymath}
The surface $F_0$ does not have a section with negative
self-intersection; however, the same description holds for its Picard
group.  The canonical class of $F_r$ is

\begin{displaymath}
K_{F_r}=  -2e -(r+2)f.
\end{displaymath}

\noindent {\bf The automorphism group of $F_r$.} The automorphism group
of $F_r$, for $r>0$, surjects onto $\P GL_2 (\C)$. The kernel is the
semidirect product of $\C^*$ with $H^0(\P^1, \O_{\P^1}(r))$ where the
former acts on the latter by multiplication. Consequently, the
dimension of the automorphism group of $F_r$ is $r+5$. $F_0$ is
isomorphic to $\P^1 \times \P^1$. The automorphism group of $\P^1
\times \P^1$ is the semidirect product of $\P GL_2(\C) \times \P GL_2
(\C)$ with $\Z / 2 \Z$, hence has dimension $6$.  \smallskip

The relation between the scrolls and the Hirzebruch surfaces is
provided by
\begin{lemma} \label{hirscr}
  The scroll $S_{k,l}$ is the image of the Hirzebruch surface
  $F_{l-k}$ under the complete linear series $|\O_{F_{l-k}}(e + lf)|$.
  For $k \not= 0, l$, the image of the curve $E$ is the unique
  rational normal k curve on the scroll. The fibers $F$ are mapped to
  lines. Irreducible curves in the class $e+(l-k)f$ map to rational
  normal $l$ curves with linear span disjoint from the linear span of
  the image of $E$.
\end{lemma}

\noindent{\bf Section classes.} During degenerations of scrolls it
will be essential to determine the limits of their hyperplane
sections. When scrolls become reducible, their hyperplane sections
remain in section classes. 

\begin{definition}
On a Hirzebruch surface $F_r$ a cohomology class of the form $e+mf$ is
called a {\bf section class}.
\end{definition}
Irreducible curves in a section class are sections of the projective
bundle. More generally, any curve in a section class consists of a section
union some fibers.  On $S_{k,l}$ the sections of degree at most $k+l$
are rational normal curves. In particular, the irreducible hyperplane
sections are rational normal $k+l$ curves.  Our description of the
cohomology ring of $F_r$ and Lemma \ref{hirscr} imply that

\begin{lemma}\label{degree}
  A curve of degree $d$ on a scroll $S_{k,l}$ that has intersection
  multiplicity one with fibers is an element of the linear series
  $|e+(d-k)f|$ on $F_{l-k}$.
\end{lemma}

\noindent{\bf Cohomology calculations.}   Since $\O_{E}(e+mf) \cong
\O_{\P^1}(m-r)$ the long exact sequence associated to the sequence
 
\begin{displaymath}
0 \rightarrow \O_{F_r} (mf) \rightarrow  \O_{F_r} (e
+ mf) \rightarrow \O_{E} (e+mf) \rightarrow 0
\end{displaymath}
implies that if $m \geq r-1$, then

\begin{displaymath}
H^1(F_r, \O_{F_r}(mf)) \rightarrow H^1(F_r, \O_{F_r}(e+mf))
\end{displaymath}
is surjective; and if $0 \leq m \leq r-1$, then

\begin{displaymath}
h^1(\O_{F_r}(e+mf))= h^1(\O_{F_r}(mf)) + h^1(\O_{\P^1}(m-r)).
\end{displaymath}
Since $F_r$ is a rational surface $H^1(F_r, \O_{F_r})=0$. Using the
exact sequence
\begin{displaymath}
0 \rightarrow \O_{F_r} (mf) \rightarrow \O_{F_r} ((m+1)f) \rightarrow
\O_{\P^1} \rightarrow 0
\end{displaymath}
we conclude that $H^1(\O_{F_r}(mf))=0$ for $m \geq 0$ by induction on
$m$.  We, thus, compute the dimensions of all the cohomology groups
for the line bundles $\O_{F_r}(e+mf)$, $m \geq 0$. 

\begin{lemma}\label{dimsectionclasses}
  The projective dimension of the linear series $|e+mf|$ on
  $S_{k,l}$, $m\geq 0$ is given by 
\begin{equation}
r(e+mf) = \max(k-l+2m+1, \ m).
\end{equation}
\end{lemma}

\noindent {\bf Remark:} The preceding discussion proves that when
$m<l-k$ the only curves in the section classes $e+mf$ consist of the
directrix $E$ union $m-k$ fibers. However, when $m \geq l-k$, the same
dimension count implies that there must be irreducible curves in the
class $e+mf$.

\begin{lemma}\label{scrolldim}
  The dimension of the locus in the Hilbert scheme whose general point
  represents a smooth scroll $S_{k,l}$ in $\P^N$ is
  \begin{displaymath}
   (k+l+2)N+2k-4-\delta_{k,l}.
  \end{displaymath}
\end{lemma}
 
\noindent {\bf Proof}: We can think of the scrolls $S_{k,l}$ as maps from the
Hirzebruch surfaces $F_{l-k}$ into projective space $\P^N$. The map is
given by $N+1$ sections of the line bundle $\O_{F_{l-k}}(e+lf)$. This
gives  $(N+1)(k+l+2)$ dimensional choices of sections. After we
projectivize and account for the automorphism group of $F_{l-k}$,
which has dimension $l-k+5+ \delta_{k,l}$, the lemma follows.  $\Box$
\smallskip

\noindent {\bf Remark:} We defined scrolls $S_{k,l}$ as surfaces in
$\P^{k+l+1}$. In case $N < k+l+1$,  Lemma \ref{scrolldim} provides
the dimension of the projections of scrolls to $\P^N$. \smallskip

Since $\P GL(N+1)$ acts transitively on
the non-degenerate scrolls $S_{k,l}$, Kleiman's theorem assures us
that if we pick general linear subspaces $ \Lambda_i \subset \P^N$ of
codimension $c_i$ such that 
\begin{displaymath}
\sum_i (c_i -2) = (k+l+2)N+2k-4-\delta_{k,l}  
\end{displaymath}
then there will be finitely many scrolls $S_{k,l}$ meeting all
$\Lambda_i$. In the rest of the paper, we address the question of
determining this number. Since we will appeal to Kleiman's theorem
\cite{kleiman:transverse} frequently, we recall it for the reader's
convenience.

\begin{theorem}[Kleiman]\label{kleiman}
Let G be an integral algebraic group scheme, $X$ an integral algebraic
scheme with a transitive $G$ action. Let $f: Y \rightarrow X$ and $g:
Z \rightarrow X$ be two maps of algebraic schemes. For each rational
element $s$ of $G$, denote by $sY$ the X-scheme given by $y \mapsto
sf(y)$.

There exists a dense open subset $U$ of $G$ such that for every
rational element in $U$, the fibered product $(sY) \times_X Z$ is
either empty or equidimensional and its dimension is the expected
dimension
\begin{displaymath}
dim(Y) + dim (Z) - dim (X).
\end{displaymath} 
Furthermore, for a dense open set this fibered product is regular.
\end{theorem}
{\bf Remark:} Although we stated the theorem in the language of
schemes, its proof holds without change for algebraic stacks.

\section{A Compactification of the Space of Scrolls}

In this section we describe a compactification of the space of
rational scrolls given by the Kontsevich moduli space of genus 0
stable maps to the Grassmannian.  \smallskip

\noindent {\bf Scrolls as curves in the Grassmannian.} To study the
geometry of scrolls it is useful to think of them as rational
curves in the Grassmannian $\G(1,N)$ of lines in $\P^N$. 
\smallskip

$S_{k,l}$ is a projective bundle over $\P^1$. The fibers of the
projection map $\pi : S_{k,l} \rightarrow \P^1$ are lines in $\P^N$.
Hence, $\pi$ induces a rational curve of Pl\"ucker degree $k+l$ in
$\G(1,N)$. More explicitly, consider the incidence correspondence

\begin{displaymath}
 \Phi = \{ (p,[L]): p \in \P^1, \ [L] \in \G(1,N), \ L \subset
  S_{k,l}, \  \pi (L)
  = p \} \subset \P^1 \times \G(1,N).  
\end{displaymath}
The image of $\Phi$ under the projection of $\P^1 \times \G(1,N)$ to
the second factor gives us the required rational curve $C \subset
\G(1,N)$.

Conversely, given an irreducible, reduced rational curve $C$ of degree
$k+l$ in $\G(1,N)$ we can construct a rational ruled surface of degree
$k+l$ in $\P^N$. Consider the projectivization of the tautological
bundle $T$ of $\G(1,N)$ over the curve $C$

\begin{displaymath}
\Psi = \{ ([L_c], p) : p \in L_c, \ c \in C \} \subset C \times \P^N
\subset \G(1,N) 
\times \P^N.  
\end{displaymath}
Projection to the second factor gives a surface $S$ of degree $k+l$ in
$\P^N$.  If the span of $S$ is $\P^{k+l+1}$, then by Proposition
\ref{classification} the surface is a rational normal scroll.  If the
span of $S$ is smaller, then $S$ is the projection of a rational
normal scroll from a linear subspace of $\P^N$. \smallskip

\noindent {\bf Non-degenerate curves.} The span of the surface has
dimension smaller than $k+l+1$ if and only if the curve is contained
in a $\G(1,r)$ for some $r < k+l+1$. We will refer to rational curves
$C\subset \G(1,N)$ which do not lie in any $\G(1,r)$ for $r < k+l+1$
as {\it non-degenerate rational curves in the Grassmannian}.
\smallskip

\noindent {\bf Non-isomorphic scrolls of the same degree.} The
correspondence between rational curves in $\G(1,N)$ and scrolls in
$\P^N$ does not yet differentiate between non-isomorphic scrolls that
have the same degree.  The automorphism group of $\G(1,N)$ does not
act transitively on non-degenerate rational curves.  The restriction
of the tautological bundle $T$ of $\G(1,N)$ to different curves can
have different splitting types. Let $\phi: \P^1 \rightarrow C$ be the
normalization of $C$. Consider the vector bundle $V= \phi^* T$ on
$\P^1$.

\begin{definition}
  We define the degree $k$ of the summand of minimal degree in the
  decomposition of $V \rightarrow \P^1$ to be the {\bf directrix
    degree} of $C$.
\end{definition} 
The directrix degree distinguishes curves associated to non-isomorphic
scrolls. Suppose $C \subset \G(1,N)$ is an irreducible, non-degenerate
curve of directrix degree $k$, then the projectivization of $V
\rightarrow \P^1$ is isomorphic to $F_{l-k}$.  The reverse
construction shows that the curve associated to the scroll $S_{k,l}$
has directrix degree $k$. We conclude that there is a natural
bijection between the set of scrolls $S_{k,l}$ in $\P^N$ and the set of
non-degenerate rational curves of degree $k+l$ and directrix degree
$k$ in $\G(1,N)$.
\smallskip

\noindent {\bf $S_{0,1}$ and $S_{1,1}$.} Unlike other scrolls, $\P^2$
and a smooth quadric $Q \subset \P^3$ have more than one scroll
structure. $\P^2$ can be given the structure of $S_{0,1}$ in a two
parameter family of ways depending on the choice of the vertex point.
The quadric surface has two distinct $S_{1,1}$ structures depending on
the choice of ruling on the quadric surface. The correspondence
between scrolls and rational curves in the Grassmannian differentiates
between these scroll structures.  \smallskip

\noindent {\bf A compactification of the space of scrolls.}
Using the preceding discussion we can compactify the space of
$S_{k,l}$ using the Kontsevich space of stable maps.  \smallskip

Let $\overline{\mathcal{S}} \subset \mbox{Hilb} (\P^N , \frac{k+l}{2}
\ x^2 + \frac{k+l+2}{2} \ x + 1)$ denote the component (with its
reduced induced structure) of the Hilbert scheme which parametrizes
rational normal scrolls. Let $\mathcal{S} \subset
\overline{\mathcal{S}}$ denote the open subscheme whose points
represent reduced, irreducible, non-degenerate scrolls.  Let
$\mathcal{C} \subset \overline{M}_{0,0} (\G(1,N), k+l)$ be the locus
in the Kontsevich moduli scheme of stable maps whose points represent
injective maps from an irreducible $\P^1$ to a non-degenerate curve in
$\G(1,N)$ of Pl\"ucker degree $k+l$.  This locus is contained in the
automorphism free locus.

\begin{theorem}\label{scrcur}
  When $k+l > 2$, there is a natural isomorphism between $\mathcal{S}$
  and $\mathcal{C}$ taking the locus of $S_{k,l}$ to maps to
   curves of directrix degree $k$.
\end{theorem}

\noindent {\bf Proof:} Projection to the second factor from the
incidence correspondence $\Phi$ induces a morphism from $\mathcal{S}$
to $\mathcal{C}$. We already observed that this morphism is a
bijection on points. Since $\mathcal{C}$ is a smooth, quasi-projective
variety \cite{fultonpan:kontsevich}, Zariski's Main Theorem implies
that this morphism is an isomorphism. $\Box$ \smallskip

{\bf Remark 1.} When $k+l \leq 2$, Theorem \ref{scrcur} is still valid
if instead of the Hilbert scheme we use the space of pointed planes
for $k+l=1$ and the space of quadric surfaces with a choice of ruling
when $k+l=2$. 
\smallskip

{\bf Remark 2.} Theorem \ref{scrcur} implies that $\mathcal{S}$ is
smooth. Note that the Hilbert scheme can be singular along subloci of
$\mathcal{S}$. For example, the Hilbert scheme of quartic scrolls is
singular along the locus of rational quartic cones---the component
corresponding to Veronese surfaces meets the component of scrolls
along that locus.

Theorem \ref{scrcur} provides us with a compactification of the space
of scrolls $S_{k,l}$. For balanced scrolls we can take the Kontsevich
moduli space of stable maps. A Zariski-open set of $\overline{M}_{0,0}
(\G(1,N), k+l)$ corresponds to maps from an irreducible $\P^1$ to a
curve of directrix degree $\lfloor \frac{k+l}{2} \rfloor$. For other
scrolls we take the scheme-theoretic closure of the locus of maps from
$\P^1$ to $\G(1,N)$ whose image has directrix degree $k$.

\section{Limits of Scrolls in One Parameter Families}

In this section we describe the limits of scrolls and their section
classes in one parameter families. 
\medskip

\noindent {\bf One parameter families of scrolls.} Let $\X \rightarrow
B$ denote a flat family of surfaces over a smooth, connected base
curve $B$. We assume that every member of the family except for the
central fiber $\X_0 \rightarrow b_0 \in B$ is a scroll $S_{k,l}$. To
simplify the statements we assume that the surface underlying $\X_0$
is still non-degenerate. This assumption can be weakened by
considering projections.

\begin{proposition}\label{centralfiber}
  The special fiber $\X_0$ is a connected surface whose irreducible
  components are scrolls $S_{k_i,l_i}$.  $\X_0$ is the image of a
  union of Hirzebruch surfaces whose dual graph is a connected tree.
  The indices $k_i, l_i$ satisfy the constraints:
\begin{equation*}
\begin{split}
&1. \ \sum_i (k_i + l_i) = k+l \\
&2. \  \sum_i k_i \leq k
\end{split}
\end{equation*}
\end{proposition}
\noindent {\bf Proof:} Since the family $\X \rightarrow B$ is flat the
central fiber $\X_0$ has to be a connected surface of degree $k+l$.
The family $\X \rightarrow B$ gives rise to a family of curves $\Y
\rightarrow B \backslash b_0$ in $\G(1,N)$, hence to a curve in
$\overline{M}_{0,0}(\G(1,N), k+l)$. Since the latter is complete we
can extend the family over $b_0$ by a stable map to $\G(1,N)$. The
projectivization of the pull back of the tautological bundle maps to
$\P^N$ giving a family that agrees with $\X$ except possibly at
$\X_0$. There is a unique scheme structure on the image that makes the
family flat. Since over a smooth base curve there is a unique way to
complete a family to a flat family (\cite{hartshorne:book} III.9.8),
this family must coincide with our original family. Hence the
underlying surface of $\X_0$ is the image of a tree of Hirzebruch
surfaces. Since $\X_0$ is non-degenerate every component has maximal
span, so by Proposition \ref{classification} is a scroll. In fact,
these scrolls must form a tree except more than 2 components might
contain a single fiber. By comparing Hilbert polynomials we can also
see that $\X_0$ does not have any embedded subschemes.

If $k=l$, relation 1 implies inequality 2.  Hence, to prove that $\sum
k_i \leq k$, we can assume that $k<l$. Then $S_{k,l}$ has a unique
degree $k$ rational curve meeting all the fibers or is a cone. The
flat limit of the degree $k$ curve is again an arithmetic genus 0
curve of degree $k$. Since meeting the fibers is a closed condition,
the limit curve has to meet all the fibers. The smallest degree curve
meeting all the fibers on $S_{k_i,l_i}$ has degree $k_i$. Hence,
$\sum_i k_i \leq k$. $\Box$  \smallskip

\noindent Proposition \ref{centralfiber} raises the question of which
unions of scrolls can be the limits of $S_{k,l}$. We now describe two
standard constructions of families of $S_{k,l}$ breaking into a
collection of $S_{k_i,l_i}$. Using these inductively we can
degenerate $S_{k,l}$ to a tree of surfaces with any $k_i,l_i$
satisfying the numerical conditions of Proposition \ref{centralfiber}.
However, we cannot smooth all such trees to an $S_{k,l}$. 
\smallskip

\noindent {\bf Example:} Cones provide the simplest
counterexample. The limit of a family of cones is a union of cones
whose vertices coincide. We can take two quadric cones meeting along a
line, but whose vertices are distinct. This surface cannot be deformed
to an $S_{0,4}$. \smallskip

\noindent {\bf Construction 1.} For any $k \geq r \geq 0$ there exists
a flat family of scrolls $S_{k,l}$ specializing to $S_{k-r,l+r}$. To
construct such a family it suffices to exhibit a flat family of vector
bundles $\O_{\P^1}(k) \oplus \O_{\P^1}(l)$ degenerating to
$\O_{\P^1}(r) \oplus \O_{\P^1}(k+l-r)$. Since $r < k$ there exists a
non-trivial injective bundle map from $\O_{\P^1}(r)$ to $\O_{\P^1}(k)$
giving rise to the extension

\begin{displaymath}
0 \rightarrow \O_{\P^1}(r) \rightarrow \O_{\P^1}(k) \oplus
\O_{\P^1}(l) \rightarrow \O_{\P^1}(k+l-r) \rightarrow 0.
\end{displaymath}
This extension gives rise to a family $E_t$ of vector bundles whose
general member is $\O_{\P^1}(k) \oplus \O_{\P^1}(l) $, but $E_0 \cong
\O_{\P^1}(r) \oplus \O_{\P^1}(k+l-r)$.  Pick the one-dimensional
subspace of $H^1 ( \P^1, \O_{\P^1}(-(k+l-r))\otimes \O_{\P^1}(r))$
which contains the extension in question. This provides us with a
family which gives $\O_{\P^1}(k) \oplus \O_{\P^1}(l) $ when $t \not=
0$ and $\O_{\P^1}(r) \oplus \O_{\P^1}(l+k-r)$ when $t=0$.

For a more geometric description of a family of $S_{k,l}$ degenerating
to $S_{k-1,l+1}$ consider a surface $S_{k,l+1}$.  When we project the
surface from a point away from the directrix, we obtain $S_{k,l}$.
However, when we project the surface from a point on the directrix we
obtain $S_{k-1,l+1}$. Now projecting $S_{k,l+1}$ from the points along
a curve that meets the directrix in isolated points, we obtain the
desired family. This construction easily generalizes to $r >1$.
\medskip

\noindent{\bf Construction 2.} There exists a family of scrolls
$S_{k,l}$ degenerating to the union of $S_{k_1,l_1}$ and $S_{k_2,l_2}$
with $k_1 + k_2 = k$.  We think of scrolls as projectivizations of
vector bundles over $\P^1$.  We choose a flat family of $\P^1$s with
smooth total space over the unit disk whose general member is smooth,
but whose central fiber has two components meeting transversely at one
point. Given a line bundle $\O_{\P^1}(k)$ on the general fiber there
is always a limit line bundle on the special curve. However, the limit
does not have to be unique. Limits differ by twists of one of the
components of the reducible fiber. We can get any splitting of $k$ on
the two components. A similar consideration applies for
$\O_{\P^1}(l)$.  Taking the desired splitting and projectivizing gives
us the desired family of Hirzebruch surfaces.  \smallskip

\noindent{\bf Remark:}  Since $\overline{\mathcal{M}}_{0,0}(\G(1,N),
k+l)$ is an irreducible, smooth Deligne-Mumford stack, every union of
scrolls whose dual graph is a connected tree can be smoothed to an
$S_{k,l}$ for some $k$ and $l$.  However, the minimal $k$ depends on
the alignment of directrices on the reducible surface (see Example
preceding Construction 1). This is the phenomenon we would like to
analyze next.  \bigskip

\noindent {\bf The limits of section classes.} Let $\X \rightarrow B$
be a flat family of scrolls subject the hypotheses in the first
paragraph of \S 4. Let $\mathcal{C} \rightarrow B$ be a flat family of
curves such that $\mathcal{C}_b \subset \X_b$ is a smooth curve in a
section class for $b \not= b_0$. We say a curve on $S_{0,l}$ is in a
section class if it is the image of a curve in a section class on
$F_l$.
 
\begin{proposition}\label{limitsection}
  The limit $\mathcal{C}_{b_0}$ restricts to a section class on each
  component of $\X_0$.
\end{proposition}

\noindent {\bf Proof: } By Proposition \ref{centralfiber} the central
fiber $\X_0$ is the union of scrolls, so it is meaningful to require
the restriction of a curve to a component to be in a section class.
Since meeting the fibers is a closed condition, $\mathcal{C}_{b_0}$
meets each fiber. To see that it does not meet the general fiber more
than once (away from the cone point of any $S_{0,l}$), consider the
one parameter family $\Y \rightarrow B$ of rational curves in
$\G(1,N)$ corresponding to our family of surfaces. Every component of
the central fiber in this family is reduced. Hence, the total space of
the family cannot be singular along an entire component. Each curve
$\mathcal{C}_b$, for $b \not= b_0$, maps isomorphically to $\Y_b$. By
Zariski's Connectedness Theorem the fibers of the map are connected
over the smooth locus. The proposition follows. $\Box$

\smallskip

\noindent {\bf Fact from Intersection Theory.} Finally, we recall a
fact from intersection theory (see \cite{fulton:intersection} chapter
12) that will be helpful in determining limits of section classes. Let
$\mathcal{C}_1$ and $\mathcal{C}_2$ be two flat families of curves
contained in a flat family of surfaces $\X \rightarrow B$ over a
smooth base curve. Assume that the general fiber of the family is
smooth and that on the general fiber $\mathcal{C}_{1,b}$ and
$\mathcal{C}_{2,b}$ are smooth curves that meet transversely at
$\gamma$ points. Let $I_0 \subset \X_0$ denote the set of isolated
points of intersection of $\mathcal{C}_{1,0}$ and $\mathcal{C}_{2,0}$
contained in the smooth locus of $\X_0$. Let $i_{p}$ denote the
intersection multiplicity at the point $p$.

\begin{lemma}\label{intersection} The intersection multiplicities
  satisfy the inequality
  $$
  \gamma \geq \sum_{p \in I_0} i_p( \mathcal{C}_{1,0},
  \mathcal{C}_{2,0}).$$
\end{lemma}

\noindent {\bf Limits of directrices.}  Proposition \ref{limitsection}
allows us to determine the limits of directrices as scrolls
degenerate. To ease the exposition first assume that the central fiber
of $\X \rightarrow B$ has two components $S_{k_i,l_i}$, $i=1,2$. For
definiteness let $k_1 + l_2 \leq k_2 + l_1$. We still assume that the
limit surface spans $\P^{k+l+1}$. Whenever we refer to the directrix of a
perfectly balanced scroll, we mean `a' directrix.

\begin{proposition}\label{direc}
  The flat limit of the directrices is a connected curve of total degree
  $k$ whose restriction to each surface is in a section class.
  Conversely, any connected curve $D$ of degree $k$, $k \leq \sum (k_i +
  l_i)/2$, whose restriction to each of the surfaces is in a section
  class is the limit of the directrices of a one parameter family of
  scrolls $S_{k,l}$.
\end{proposition}

\noindent {\bf Proof:} The first assertion is a restatement of
Proposition \ref{limitsection}.  To prove the second assertion we will
explicitly construct the desired families of $S_{k,l}$.  There are two
cases depending on whether $D$ contains a multiple of the line $L$
joining the two surfaces. If $D$ does not contain $L$, then $D$ must
consist of a section (possibly with some fibers) in each surface
meeting $L$ at the same point.  \smallskip

When $S_{k,l}$ degenerates to $S_{k-1,l+1}$, the limit of the
directrices of $S_{k,l}$ is the union of the directix of $S_{k-l,l+1}$
and a fiber since these are the only connected degree $k$ curves in a
section class on $S_{k-1,l+1}$.  Since the projective linear group
acts transitively on the fibers of a scroll, using Construction 1 we
conclude that the union of the directrix and a fiber on $S_{k-1,l+1}$
can be smoothed to the directrix of an $S_{k,l}$. Inductively, we
conclude the analogous statement for $S_{k-r,l+r}$. We can, therefore,
assume that $D$ consists of 2 sections meeting $L$ at the same point.
\smallskip

  If $k=k_1 + k_2$, then the limit of the directrices
must be the union of the directrices of the limit surfaces. A section
class on a scroll $S_{k_i,l_i}$ has degree at least $k_i$.  Since the
total degree is $k = \sum_i k_i$, the curve must break exactly into
degree $k_i$ section classes. A section class of degree $k_i$ on
$S_{k_i,l_i}$ is the directrix. Using Construction 2 and the fact that
$\P GL(k+l+2)$ acts transitively on the pairs of surfaces
$S_{k_i,l_i}$ that span $\P^{k+l+1}$, meet along a line and whose
directrices meet along their common line, we conclude that such pairs
can be smoothed to $S_{k,l}$ so that the union of their directrices
deforms to the directrix of $S_{k,l}$. 
\smallskip

{\bf The case $k_1 + k_2 < k$.} Since we are assuming that $k_1+l_2
\leq k_2 + l_1$, $D$ must consist of the directrix in $S_{k_1,l_1}$
and a section of degree $k-k_1$ in $S_{k_2,l_2}$. Set
$j=l+k-k_1-l_1-k_2$. Consider a family of $S_{k, l+ k-k_1-k_2}$ over a
small disk in the complex plane degenerating to the union of
$S_{k_1,l_1}$ and $S_{k-k_1, m}$ as described in Construction 2. By
our discussion the directrices must specialize to the union of the
directrices. Now choose $k-k_1-k_2$ general sections of the family
that all pass through general points of $S_{k-k_1, j}$. (Observe that
the total space of the family we exhibited in Construction 2 is
generically smooth on every component of the special fiber, so we can
select such sections.) Projecting the family along these sections
gives a family of $S_{k,l}$ having the desired numerical properties.

We need to verify that we can get all sections $C$ of degree $k-k_1$
on $S_{k_2,l_2}$ as the projection of the directrix of $S_{k-k_1, j}$
from suitable points on the surface.  On $S_{k_2,l_2}$ blow up
$k-k_1-k_2$ general points $p_i$ on $C$. Let $\Xi$ be the sum of the
exceptional divisors.  The linear series $e + (k-k_1-k_2)f - \Xi$ maps
the blow up to projective space as $S_{k-k_1, j}$. This map contracts
the fibers passing through $p_i$ and maps $C$ to the directrix.
Projecting $S_{k-k_1, j}$ from the points corresponding to the image
of the contracted fibers projects $S_{k-k_1, j}$ to $S_{k_2,l_2}$ and
the directrix onto $C$. This concludes the construction.  \medskip

Now we treat the case when $D$ contains a multiple of the common line
$L$. As in the previous case we can assume that $D$ consists of two
sections and $L$ with multiplicity $m$. $D$ must contain the directrix
in $S_{k_1,l_1}$ and a section of degree $k-k_1-m$ in $S_{k_2,l_2}$.
Using the previous case inductively, we can find a family of scrolls
$S_{k+m,l+m}$ degenerating to a chain of $S_{k_1,l_1}$, $m$ quadric
surfaces and $S_{k_2,l_2}$ such that their directrices specialize to
the directrix of $S_{k_1,l_1}$, a section of degree $k-k_1-m$ in
$S_{k_2,l_2}$ and a chain of conics on the quadric surfaces connecting
these two curves. Choose $2m$ general sections which specialize to a
pair of points on a fiber on each quadric surface. Projecting the
family from those sections gives the desired family. $\Box$
\smallskip

\noindent {\bf Remark:} We can inductively extend Proposition
\ref{direc} to the case when the central fiber contains more than two 
components. The following theorem summarizes the conclusion:

\begin{theorem}
  Suppose a one parameter family of scrolls $S_{k,l}$ specializes to
  the union $\bigcup_{i=1}^r S_{k_i,l_i}$. Then the limit of the
  directrices is a connected curve of degree $k$ whose restriction to
  each surface is in a section class. Conversely, given any connected
  curve $C$ of degree $k \leq \sum_i (k_i +l_i)/2$ whose restriction
  to each component is in a section class, there exists a one
  parameter family of $S_{k,l}$ specializing to the reducible surface
  such that the limit of the directrices is $C$.
\end{theorem}

\noindent {\bf Limits of other section classes.} 
When a family $\X \rightarrow B$ of scrolls $S_{k,l}$ specializes to a
union of two scrolls $S_{k_1,l_1} \bigcup S_{k_2, l_2}$, then the flat
limit of curves in a section class $e+mf$, $m \geq l-k$, are connected
curves of total degree $m+k$ that restrict to section classes. Suppose
that the total space of the family is smooth. The curves give a line
bundle $L$ over $\X \backslash \X_0$. We can always extend this line
bundle to the entire family. However, when the central fiber is
reducible, this extension is not unique. Twisting by the components of
the central fiber give different extensions.

For concreteness, suppose the line bundle $L$ is the pull-back of
$\O_{\P^N}(1)$ to $\X \backslash \X_0$. Let $L_0$ be the line bundle
over $\X_0$ arising from any extension of $L$. If the restriction of
$L_0$ to each component is effective, then the degree of $L_0
|_{S_{k_1,l_1}}$ ranges between $k_1$ and $k+l - k_2$. One component
of the limit curves corresponds to hyperplanes not containing either
of the components of the limit scroll. The other limits correspond to
hyperplane sections by hyperplanes containing one of the scrolls and
tangent to a certain order to the other one along their common line.
\smallskip

\noindent {\bf Example:} When a family of scrolls $S_{2,4}$
specialize to $S_{1,1} \cup S_{1,3}$, a possible limit of the
hyperplane sections restricts to a degree $5$ curve on $S_{1,1}$ and
the directrix on $S_{1,3}$. The dimension of such curves is $8$.
However, the dimension of the hyperplane sections was only 7.
Consequently, not all quintics can arise as the limit of hyperplane
sections of our family of scrolls.  Lemma \ref{intersection} provides
the answer. On the smooth surfaces the directrices have intersection
number 2 with the hyperplane sections. The limit of the directrices is
the union of the directrix $D$ of $S_{1,3}$ and the directrix $L$ on
$S_{1,1}$ meeting $D$. A general quintic meeting $D$ meets $L$
transversely in 3 other points. We conclude that a quintic can be part
of a limit of the hyperplane section only if it has contact of order 2
with $L$ at their point of intersection on the common fiber.

The following lemma, which is an immediate consequence of Lemma
\ref{intersection} and Proposition \ref{limitsection}, summarizes the
general situation.

\begin{lemma}
  Suppose a one parameter family of scrolls $S_{k,l}$ specializes to a
  union $\bigcup_{i=1}^r S_{k_i,l_i}$. Then the limit of curves in a
  section class $e+mf$ specialize to connected curves of degree $m+k$.
  Their restrictions to each surface lie in a section class.
  Furthermore, the sum of the intersection multiplicities of a limit
  curve with the limit of the directrices at isolated points of their
  intersection on the smooth locus of the surface cannot exceed
  $m+k-l$.
\end{lemma}

\section{Examples}

\noindent {\bf Example A. Counting cubic scrolls in $\P^4$.} Since by Lemma
\ref{scrolldim} there is an 18 dimensional family of cubic scrolls
$S_{1,2}$ in $\P^4$, there are finitely many containing $9-n$ general
points and meeting $2n$ general lines.  \smallskip

\begin{figure}[htbp]
\begin{center}
\epsfig{figure=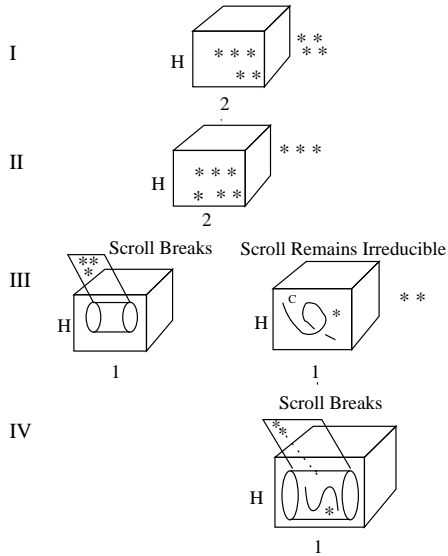}
\end{center}
\caption{Example A1. Cubic scrolls containing 9 general points in $\P^4$.}
\label{Figure 2}
\end{figure}

\noindent {\bf A1. Cubic scrolls containing 9 general points.}  We
specialize the nine points one by one to a fixed hyperplane $H$ in
$\P^4$ (see Figure \ref{Figure 2}). We can take $H$ to be the span of
four of the points.  \smallskip

\noindent {\bf Step I.} Specialize a fifth point $p_5$ to a general
point of $H$. There are no reducible scrolls at this stage containing
all the points. Any reducible scroll would be the union of a quadric
and a plane meeting along a line.  Since no six of the points lie in a
$\P^3$, the quadric could contain at most $5$ of the points. However,
the remaining 4 points do not lie on a plane.  \smallskip

\noindent {\bf Step II.}  Specialize a sixth point $p_6$ to $H$. Now
there is a reducible solution: the plane $P$ spanned by the 3 points
outside $H$ and the unique quadric in $H$ containing $H \cap P$ and
the six points in $H$. However, at this stage there might still be
irreducible solutions. Their hyperplane section in $H$ must be the
unique twisted cubic $C$ that contains the 6 points in $H$.
\smallskip

\noindent {\bf Step III.}  Specialize a seventh point
$p_7$ to a general point of $H$.  Bezout's theorem forces the scrolls
to break into a union of a quadric surface and a plane. The quadric
$Q$ must contain the twisted cubic $C$ and $p_7$. Since the plane and
the quadric meet in a line, $Q$ must also contain the point of
intersection $q$ of $H$ with the line spanned by the points outside
$H$.  This determines the quadric uniquely. The plane is also
determined because it must contain the line in the quadric through $q$
which meets $C$ only once: recall by Lemma \ref{limitsection} the
curve $C$ is in a section class.

Later we will check that both of the solutions occur with
multiplicity 1. This will prove that there are 2 cubic scrolls
containing 9 points in $\P^4$.  
\smallskip

\noindent {\bf A2. Cubic scrolls in $\P^4$ containing 6 points and meeting
  6 lines.} Here we sketch the degenerations required to see that
there are 1140 cubic scrolls in $\P^4$ containing 6 points and meeting
6 lines until we reduce the problem to a straightforward problem about
quadrics and planes (see Figure \ref{Figure 3}). \smallskip

\noindent {\bf Step I.} We specialize 5 points and a line $l_1$ to a
fixed hyperplane $H$. This is the first stage where reducible
solutions exist. There can be a quadric $Q$ contained in $H$ and a
plane $P$ outside $H$ meeting $Q$ in a line.  There are 4
possibilities: of the 5 lines outside $H$, 4, 3, 2, or 1 of them can
meet $Q$ and 1, 2, 3, or 4 remaining lines, respectively, can meet
$P$.  

If 4 of the lines meet $Q$, then $Q$ is determined uniquely. $P$ can
be any of the 4 planes meeting the remaining line, containing the
remaining point and meeting $Q$ in a line. Since $l_1$ meets $Q$ in 2
points, each of the solutions count twice for the choice of point.
Finally, we have a factor of 5 for the choice of which 4 lines among
the 5 meet the quadric $Q$. We express this as $5 \times 2 \times
4=40$ where the first multiple is the combinatorial muliplicity for
the choice of lines, the second multiple is for the choice of point
and the last number is the number of surfaces satisfying the incidence
conditions.  The analysis of the other three cases is similar.

At this stage some scrolls can remain irreducible. Their hyperplane
section $C$ in $H$ is then a twisted cubic containing the 5 points and
meeting $l_1$.  \smallskip

{\noindent {\bf Step II.}  We specialize another line $l_2$ to $H$.
  There are new reducible solutions.

{\bf Case i.} There can be a plane $P$ in $H$ and a quadric $Q$
meeting it along a line. $P$ must be one of the 10 planes spanned by 3
of the points in $H$.  Finally, $Q$ must contain the other two points
in $H$, meet $P$ in a line, and meet the lines and the point $p$ lying
outside $H$. Further specialization shows that there are 6 such
quadrics. Hence we get 60 solutions.

{\bf Case ii.} There can be a quadric $Q$ in $H$ and a plane $P$
outside meeting $Q$ in a line. There are four possibilities: 1, 2, 3
or 4 of the lines can meet $P$.  Let us analyze the case when 3 lines
meet $Q$. The limit hyperplane section $C$ contained in $Q$ must meet
5 points and $l_1$. Among the the quadrics containing 8 points only a
finite number contains such a twisted cubic. To determine the number
we specialize the conditions on the cubic curve.

Take a general 2 plane $\Pi$ and specialize 3 of the points and $l_1$
to $\Pi$. By Bezout's theorem $C$ has to break into a conic and a
line. The line can be any of the three lines containing two of the
points in $\Pi$ or it can be the line joining the two points outside
of $\Pi$. Once we require $Q$ to contain any of the lines, it is
uniquely determined.  When there is a line $l$ in $\Pi$, the limit of
$C$ meets $l_1$ only in $l \cap l_1$. When the line $l$ is outside
$\Pi$, the limit of $C$ has a conic in $\Pi$ which must meet $l_1$ in
2 points, hence we get a multiplicity of 2. The other cases are
analogous.

\begin{figure}[htbp]
\begin{center}
\epsfig{figure=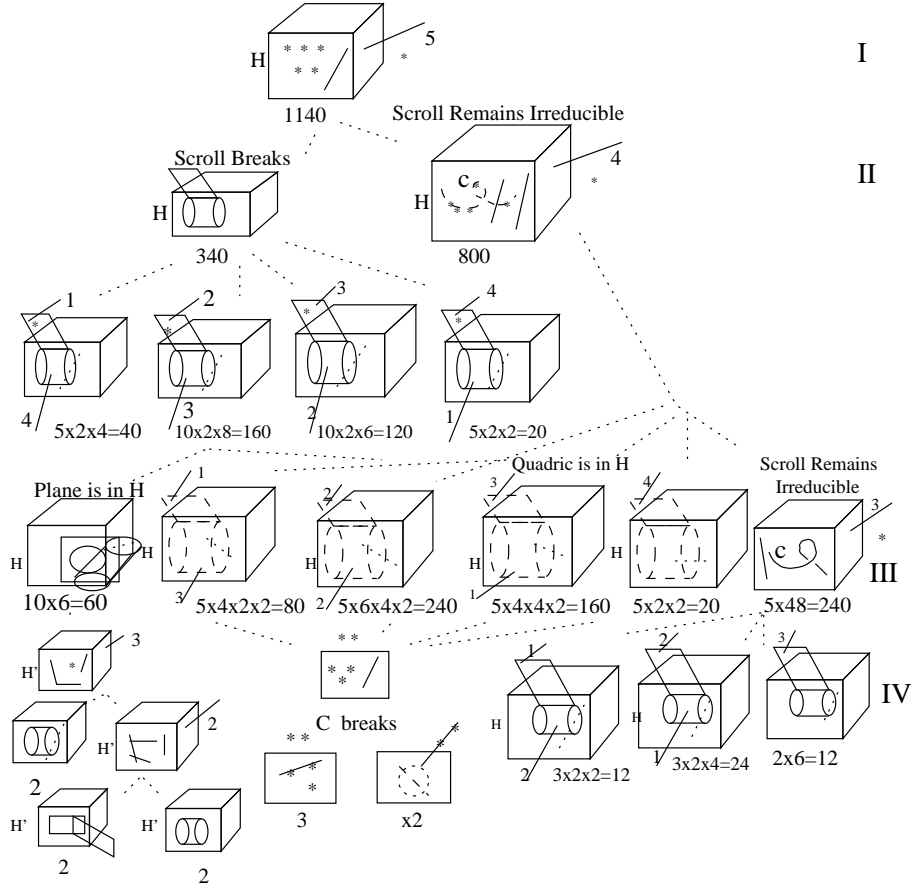}
\end{center}
\caption{Example A2. Cubic scrolls in $\P^4$ containing 6 points and
  meeting 6 lines.}
\label{Figure 3}
\end{figure}

By Lemma \ref{limitsection} the curve $C$ must be in a section class
in $Q$, hence should meet the fibers only once. When counting quadric
surfaces, one has to be careful to distinguish between the rulings.

{\bf Step III.} Finally, there can still be irreducible scrolls. In
that case their hyperplane section must be one of the $5$ twisted
cubics in $H$ containing 5 points and meeting $l_1$ and $l_2$ (see \S
2.3 of \cite{vakil:rationalelliptic}). The analysis is similar to the
previous cases.  \smallskip

\noindent {\bf Example B. Counting quadric surfaces in $\P^4$.} The
degeneration method allows us to count different types of scrolls.  We
illustrate this by counting quadric surfaces and quadric cones in
$\P^4$ (see Figure \ref{Figure 4}).  \smallskip

\noindent {\bf B1. Quadric surfaces in $\P^4$ containing 3 points and
  meeting 7 lines.} {\bf Step I.}  Specialize the three points
$p_1,p_2,p_3$ and a line $l_1$ to the hyperplane $H$. At this stage
there is a unique quadric contained in $H$ satisfying all the
incidence conditions. It counts with multiplicity 2 for the choice of
intersection point with $l_1$.

If a quadric is not contained in $H$, its hyperplane section must lie
in the plane $\Pi$ spanned by the points $p_1, p_2, p_3$ and must meet
$l_1$ at $l_1 \cap \Pi$.

{\bf Step II.} We specialize $l_2$. The quadric can lie in $H$. If
not, the hyperplane section must be the unique conic containing $p_i$
and $l_j \cap \Pi$. Specializing a third line $l_3$ forces the
quadrics to either become reducible or to lie in $H$. We obtain a
total of 9 quadric surfaces containing 3 points and meeting 7
lines. \smallskip

\begin{figure}[htbp]
\begin{center}
\epsfig{figure=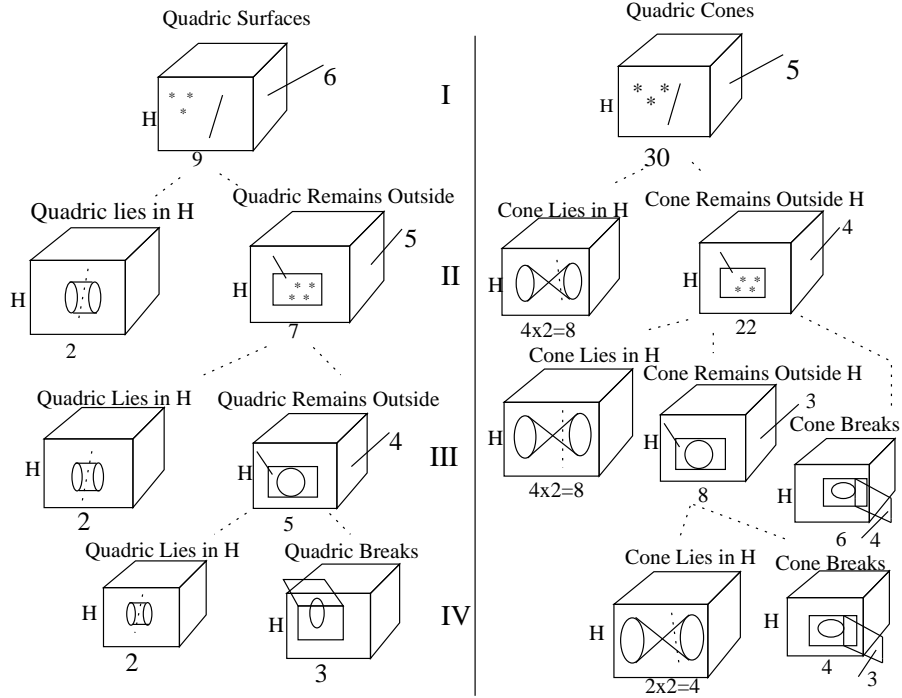}
\end{center}
\caption{Example B. Counting quadric surfaces containing 3 points and
  meeting 7 lines and quadric cones containing 3 points and meeting 6
  lines.}
\label{Figure 4}
\end{figure}

\noindent {\bf B2. Quadric cones in $\P^4$ containing 3 points and meeting 6
  lines.} We compare the case of quadric cones to the case of quadric
surfaces.

{\bf Step I.} Specialize the points $p_1,p_2,p_3$ and the line $l_1$
to $H$. The cone can lie in $H$. There are 4 quadric cones in $\P^3$
containing 8 general points. Each solution counts with multiplicity 2
for the choice of intersection point with $l_1$.

{\bf Step II.} If a cone does not lie in $H$, its hyperplane section
in $H$ must lie in the plane $\Pi$ spanned by $p_i$, so it must meet
$p_i$ and $q=l_1 \cap \Pi$. We specialize another line $l_2$ to 
$H$. The cone can lie in $H$. Again there are 4 solutions each counted
with multiplicity 2.

There can also be reducible solutions: the union of $\Pi$ and one of
the three planes that meet $\Pi$ in a line and meet the three lines
not contained in $H$. This case is delicate. The limit of the
hyperplane sections is a conic containing $p_i$ and $q$. The two
planes are images of Hirzebruch surfaces $F_1$ whose directrices are
contracted.  The image conic is in the class $e+2f$ on the
$F_1$ it lies in. It also meets the directrix of the other $F_1$.
Hence, the limit of the hyperplane section has to contain $p_i$ and
$q$ and be tangent to the line common to the planes at the limit of
the vertices of the nearby cones. There are two conics containing 4
points and tangent to a line in $\P^2$. The two points of tangency
give us the possible limiting positions of the vertices of our
original family of cones. So each of the pairs of planes can be the
limit of cones in two ways depending on the choice of the vertex
point. We thus get a count of 6.

{\bf Step III.} If the cone is neither reducible nor contained in $H$,
then the hyperplane section in $H$ is determined. Specializing a third
line $l_3$ forces the cones to break or to lie in $H$. The
calculations  are analogous to the previous case. We obtain
a total of 30 quadric cones containing 3 points and meeting 6 lines.
\smallskip

\noindent {\bf Example C. Counting scrolls of degree $n$ in $\P^{n+1}$
  containing $n+5$ points and meeting an $n-3$ plane.} We give a final
example to illustrate the types of recursive formulae one might hope
to obtain from our method of counting. Observe that by Lemma
\ref{scrolldim} there will be finitely many scrolls $S_{\lfloor
  \frac{n}{2} \rfloor, \lfloor \frac{n+1}{2} \rfloor}$ in $\P^{n+1}$
containing $n+5$ points and meeting an $n-3$ plane $P$. Let us denote
this number by $S(n)$.

\begin{figure}[htbp]
\begin{center}
\epsfig{figure=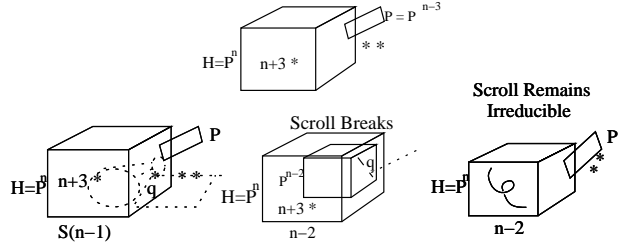}
\end{center}
\caption{Counting degree $n$ scrolls in $\P^{n+1}$ containing $n+5$ points
  and meeting an $n-3$ plane.}
\label{Figure 5}
\end{figure}

{\bf Step I.} We specialize the points to a hyperplane $H$. The first
irreducible solution occurs when $n+3$ of the points lie in $H$: a
plane outside $H$ and a scroll of degree $n-1$ in $H$ meeting the
plane in a line. The degree $n-1$ scroll has to meet $q$, the point of
intersection of $H$ with the span of the two points outside $H$. If
the scroll meets $P$, then we reduce to the same problem in degree one
less, so the number is $S(n-1)$.

If the plane meets $P$, then the scroll must contain a line $l$ in
$P'= H \cap P$ containing $q$. We have to count scrolls of degree
$n-1$ in $\P^n$ containing a line through a point in a $\P^{n-2}$ and
containing an additional $n+3$ points in general position. An easy
specialization shows that there are $n-2$ such scrolls.

If the scrolls remain irreducible, then their hyperplane section $C$
in $H$ is determined. When we specialize the $n-3$ plane $P$ to $H$
the scroll breaks into a degree $n-1$ scroll union a plane. We are
reduced to counting degree $n-1$ scrolls containing $C$, a point and
meeting $P$. It is easy to see that there are $n-2$ such scrolls.
Solving the recursion we conclude that there are $(n-1)(n-2)$ degree
$n$ scrolls in $\P^{n+1}$ containing $n+5$ points and meeting an $n-3$
plane.

We will now justify the calculations made above by making the
necessary dimension counts and multiplicity calculations. In a table
at the end of the paper we will provide some other characteristic
numbers of surfaces.

\section{Degenerations Set Theoretically}

In this section we describe the set theoretical limits of surface
degenerations under the assumption that the surfaces and their
successive hyperplane sections remain non-degenerate.  We will compute
the dimension of images of maps from trees of Hirzebruch surfaces to
$\P^N$ and determine the codimension one loci in the Hilbert scheme of
scrolls when we require one of their points to lie in a fixed
hyperplane. The calculation in this section will be purely set
theoretic. We collect our notation here for the reader's convenience.
\smallskip

\noindent {\bf Notation:} Let $H$ and $\Pi$ be general hyperplanes in
$\P^N$. In our algorithm we will specialize a linear space meeting
either a surface or a curve on a surface to a hyperplane. We will
refer to the hyperplane by $H$ when we specialize a condition on a
surface and by $\Pi$ when we specialize a condition on a curve.

The dimension of a linear space will be indicated by a subscript. We
will often omit the dimension from the notation.

Let $\{ \Delta_{a_i}^i \}_{i=1}^I$, $\{ \Sigma_{a_i'}^i \}_{i=1}^{I'}$
and $\{ \tilde{\Lambda}_{\tilde{b}_i}^i \subset \Lambda_{b_i}^i
\}_{i=1}^Y$ be three collections of general linear subspaces of
$\P^N$. These will be the linear spaces that we have not yet
specialized to the hyperplane and meet the surface, meet the marked
curve and meet and contain a fiber, respectively.

Let $\{ \Gamma_{c_j}^j \}_{j=1}^{J_H} $ and $\{ \Omega_{d_j}^j
\}_{j=1}^{J_{\Pi}}$ be collections of general linear subspaces of $H$
and $\Pi$, respectively.  These will be the linear spaces that we have
already specialized to the hyperplane.

$\Delta(0), \cdots, \Delta(M)$ will denote a partition of the $I$ linear
spaces $\Delta_i$ into $M+1$ parts. We will use analogous notation for
the other linear spaces.   
\smallskip

\noindent {\bf Spaces of maps.} We now define a sequence of spaces of
maps from Hirzebruch surfaces to $\P^N$. Intuitively
they will correspond to scrolls with two marked curves on them. The
scrolls will meet certain linear spaces. The curves will meet some
others. In addition some fibers of the scroll will be required to lie
in a linear space and meet another linear space.  \smallskip

\noindent {\bf Definition of $MS_H$.} Let

\begin{equation*}
\begin{split}
  MS_H &(\P^N;\ k,l;\ C(k+l),D; \{\lambda_i \}_{i=1}^Y, \{
  q_i
  \}_{i=1}^I , \{ p_j \in  C \}_{j=1}^{J_H}, \pi :  
  \pi^{-1} (H) = C; \\ & ( \pi(\lambda_i) \subset \Lambda^i,
  \pi(\lambda_i) \cap \tilde{\Lambda}^i \not= \emptyset)_{i=1}^Y, (
  \pi(q_i) \in \Delta^i )_{i=1}^I , ( \pi(p_j) \in \Gamma^j
  )_{j=1}^{J_H})
\end{split}
\end{equation*}
be the set of maps, up to isomorphism, from a Hirzebruch surface
$F_{l-k}$, with two marked sections $(C,D)$, $Y$ marked fibers and $I$
marked points, to $\P^N$ whose image is an $S_{k,l}$ not contained in
$H$ such that

1. $\pi(C) = \pi(S) \cap H$.

2. $D$ is the directrix (or if $k=l$, a choice of directrix).

3. The images of the marked fibers $\lambda_i$ are contained in the
linear spaces $\Lambda^i$ and meet the linear spaces
$\tilde{\Lambda}^i$.

4. $p_i$ are marked points on the surface and $q_j$ are marked points
on the curve $C$. We assume that they are distinct points whose images
lie in the linear spaces $\Delta^i$ and $\Gamma^j$, respectively.
\smallskip

\noindent {\bf Definition of $MS_{\Pi}$}. Similarly let   

\begin{equation*}
\begin{split}
  &MS_{\Pi} (\P^N;\ k,l;\ C(d),D;\{\lambda_i\}_{i=1}^Y,\{ q_i
  \}_{i=1}^I , \{ q_i' \in C \}_{i=1}^{I'}, \{ p_j \in C
  \}_{j=1}^{J_{\Pi}}; \pi :
  (\pi(\lambda_i) \\
  & \subset \Lambda^i, \pi(\lambda_i) \cap \tilde{\Lambda}^i \not=
  \emptyset)_{i=1}^Y, ( \pi(q_i) \in \Delta^i )_{i=1}^I , ( \pi(q_i')
  \in \Sigma^i )_{i=1}^{I'}, (\pi( p_j) \in \Omega^j )_{j=1}^{J_{\Pi}}
  )
\end{split}
\end{equation*}
be the set of maps, up to isomorphism, from a Hirzebruch surface
$F_{l-k}$, with two marked sections $(C,D)$, $Y$ marked fibers and $I$
marked points, to $\P^N$ whose image is an $S_{k,l}$ such that

1. $\pi(C)$ has degree $d$.

2. $D$ is the directrix (or if $k=l$, a choice of directrix).

3. The images of the marked fibers $\lambda_i$ are contained in the
linear spaces $\Lambda^i$ and meet the linear spaces
$\tilde{\Lambda}^i$.

4. $q_i$ are marked points on the surface, $q_j$ are marked points on
the curve $C$ whose images are not contained in $\Pi$ and $p_j$ are
marked points on the curve $C$ that map to $\Pi$. We assume that they
are distinct points whose images lie in the linear spaces $\Delta^i$
and $\Sigma^i$ and $\Omega^j$, respectively.  \smallskip

We can compactify both $MS_H$ and $MS_{\Pi}$ in a manner analogous to
\S 3. For concreteness we explain the construction for $MS_H$. The
case of $MS_{\Pi}$ is identical.  Let $\F(0,0,1;N)$ denote the variety
of two-pointed lines in $\P^N$. A curve in $\F(0,0,1;N)$ is determined
by three degrees, the degrees $(d_0,d_1)$ of the two projections
$\gamma_0, \gamma_1$ to $\P^N$ and the degree $d_2$ of the projection
$\gamma_2$ to $\G(1,N)$.  Given a map $\pi$ in $MS$ from a surface $S$
with sections $C$ and $D$, we get a stable map to $\F(0,0,1;N)$ by
sending $(S,C,D,\pi)$ to the map $\sigma$ from $C$ given by $\sigma(p)
= (\pi(p), \pi(D \cap F_p) , [\pi(F_p)])$ where $F_p$ is the fiber
through $p$.  We have to enhance this correspondence to mark the
points on the surface. To do that we simply take the $I$-th fold fiber
product of the universal family $\mathcal{U}$ over the stack
$\overline{\mathcal{M}}_{0,J_H + Y}(\F(0,0,1;N), (d_0,d_1,d_2))$. We
will denote the closure of $MS_H$ and $MS_{\Pi}$ in these stacks as
$\overline{\mathcal{MS}}_H$ and $\overline{\mathcal{MS}}_{\Pi}$,
respectively. When we do not want to distinguish between them we will
use the notation $\overline{\mathcal{M}}$. \smallskip 

\noindent {\bf Definition of the Divisors $D_H$, $D_{\Pi}$.} In
$\overline{\mathcal{MS}}_{H}$ requiring the stable map to map $p_I$
into $H$ defines a Cartier divisor. We will denote this divisor by
$\mathcal{D}_H$. Similarly, $\overline{\mathcal{MS}}_{\Pi}$ has a
Cartier divisor $\mathcal{D}_{\Pi}$ defined by requiring the image of
$q_{I'}'$ to lie in $\Pi$. Colloquially, the surfaces in the divisors
are the surfaces we see after we specialize one of the points to a
hyperplane. \smallskip

\noindent {\bf Definitions of $X_H$ and $X_{\Pi}$.} Let
\begin{displaymath}
X_H (\P^N; (k_i,l_i;C(d_i),D(e_i), \lambda(i), q(i), p(i) \in
C(i))_{i=0,1}; \pi)    
\end{displaymath}
be the set of maps, up to isomorphism, from a pair of Hirzebruch
surfaces $F_{l_i - k_i}$ meeting transversely along a fiber with the
usual markings such that

1. $S_0 := \pi(F_{l_0 - k_0}) \subset H$ is a scroll $S_{k_0,l_0}$
contained in $H$.
   
2. $S_1 := \pi(F_{l_1 - k_1})$ is a scroll $S_{k_1,l_1}$ not contained
in $H$ and which meets $H$ transversally along the line that joins it
to $S_0$.
  
3. $C(d_i)$ is a section of degree $d_i$ and $D(e_i)$ is a section of
degree $e_i$ on $S_i$ such that $C(0) \cup C(1)$ and $D(0) \cup D(1)$
form a connected curve.
   
4. $\lambda(i)$, $q(i)$ and $p(i)$ is a partition of the marked fibers
and points to the two components and they satisfy the same incidence
and containment relations as in the definition of $MS_H$.

The definition of $X_{\Pi}$ is similar, but with $H$ replaced by $\Pi$
and the names of the linear spaces and points modified as in the
definition of $M_{\Pi}$. Denote the corresponding stacks in
$\overline{\mathcal{MS}}_H$ and $\overline{\mathcal{MS}}_{\Pi}$ by
$\mathcal{X}_H$ and $\mathcal{X}_{\Pi}$, respectively. \medskip

\noindent {\bf Definitions of $W_H$ and $W_{\Pi}$.}
Let  
\begin{displaymath}
W_H (\P^N; (k_i,l_i;C(d_i),D(e_i), \lambda(i), q(i), p(i) \in
C(d_i))_{i=0}^M; \pi)    
\end{displaymath}
be the set of maps, up to isomorphism, from the union of $M+1$
Hirzebruch surfaces with the usual markings to $\P^N$ such that
 
1. All the components $S_i$ for $i > 0$ are attached along distinct
fibers to a central component $S_0$,

2. $\pi(S_0)$ is a cone not contained in $H$,

3. $C(0)$ and $D(0)$ both contain the directrix of $S_0$,

4. $\pi (S_i) \subset H$ for all other $i$,

5. The fibers and marked points are distributed according to the
partition and meet and lie in their designated linear spaces as in the
definition of $MS_H$.

The definition of $W_{\Pi}$ is obtained by making the usual
modifications. We will denote the corresponding stacks by
$\mathcal{W}_H$ and $\mathcal{W}_{\Pi}$.  \smallskip

\noindent {\bf Subschemes of $\mbox{Hilb}(\P^N,S^{k+l})$.}  Let
$\mbox{Hilb}(\P^N,S^{k+l})$ denote the component of rational normal
scrolls of degree $k+l$ in the Hilbert scheme. Since the image
surfaces in the spaces of maps we defined have generically the same
Hilbert polynomial, they all map to $\mbox{Hilb}(\P^N,S^{k+l})$.
\smallskip

\noindent {\bf Dimension counts.} With this preparation we can begin
the dimension counts.

\begin{proposition}\label{usekl}
  Let $A$ be a reduced, irreducible substack of $\overline{\mathcal{M}}$ and
  let $p$ be any of the labeled points. Then there exists a
  Zariski-open subset $U$ of the dual projective space $\P^{N*}$ such
  that for all hyperplanes $[H] \in U$ the intersection $A \cap \{ \pi
  (p) \in H \}$ is either empty or reduced of dimension $\dim A -1$.
\end{proposition}

\noindent {\bf Proof:} $\P^N$ is a homogeneous space under the action
of the group $\P GL(N+1)$. The proposition follows from Theorem
\ref{kleiman}. In the notation of the theorem take $f: H \rightarrow
\P^N$ be the immersion of a hyperplane. Let $g : A \rightarrow
\P^N$ be evaluation morphism at $p$. $\Box$

\begin{proposition}\label{linekl}
  Let $A$ be a reduced, irreducible substack of $\overline{\mathcal{M}}$ and
  let $\lambda$ be any of the marked fibers. Then there exists a
  Zariski-open subset $U$ of the dual projective space $\P^{N*}$ such
  that for all hyperplanes $[H] \in U$ the intersection $A \cap \{ \pi
  (\lambda) \in H \}$ is either empty or reduced of dimension $\dim A
  -2$.
\end{proposition}

\noindent {\bf Proof:} Consider the Grassmannian $\G(1,N)$ under
the group action $\P GL(N+1)$. Let $f: \Sigma_{1,1}(H) \rightarrow
\G(1,N)$ be the immersion of the Schubert cycle of lines contained in $H$.
Let $g: A \rightarrow \G(1,N)$ be the evaluation morphism. The
proposition follows from Theorem \ref{kleiman} and the fact that the
cycle $\Sigma_{1,1}(H)$ has codimension 2 in the Grassmannian. $\Box$
\medskip

\noindent {\bf Remark.} In the notation of the Proposition
\ref{usekl} or \ref{linekl}, we can further deduce that if $B$ is a
proper, closed substack of $A$ then every component of $B \cap \{ \pi
(p) \in H \}$ is a proper closed substack of a component of $A \cap
\{ \pi (p) \in H \}$ by using the dimension statement in Proposition
\ref{usekl} or \ref{linekl} for every component of $B$.
\medskip

\noindent {\bf Caution and Convention:} For the rest of this section,
{\it a directrix of a perfectly balanced scroll will be considered an
  ordinary section class.} The directrix class of a perfectly balanced
scroll behaves differently than the directrix class of other scrolls
dimension theoretically. Instead of noting this exception in every theorem
{\it we declare that a perfectly balanced scroll does not have any
  directrices.}  \medskip

\noindent {\bf The dimension of the building blocks.} We now compute
the dimension of the locus of maps, up to isomorphism, from an
irreducible Hirzebruch surface $F_{l-k}$ with two marked irreducible
sections $C$ and $D$ to $\P^N$ such that
 
(i) the image of $F_{l-k}$ is a scroll $S_{k,l}$ that has contact of
order $m_i$ along fibers $\lambda_i$, $1\leq i \leq m$, with a fixed
hyperplane H,

(ii) the marked curves $C,D$ have degrees $d$ and $e$, respectively. In
case these curves are distinct we assume that they have contact of
order $n_i$ with each other along distinct fibers $\lambda_1, \cdots,
\lambda_k, \lambda_{m+1}, \cdots, \lambda_{m+n}$ ($k \leq m$).

\begin{proposition} \label{nastydim}
If $\sum_{i=1}^m m_i > k$, then this dimension is
$$N(k+l+2)+k-l-5-\delta_{k,l} + m + \alpha,\  \mbox{where}$$
 
(i) $\alpha=0$ if $C$ and $D$ lie in $H$;

(ii) $\alpha= 2d -k-l+1$ if $C$ and $D$ coincide, but do not lie in $H$;

(iii) $\alpha= 2e -k-l+1+ n - \sum_{i=1}^{n+m} n_i $ if $C$ lies in $H$,
but $D$ does not;

(iv) $\alpha = 2d+2e-2k-2l+2 +n- \sum_{i=1}^{n+m} n_i $ if $C$ and
$D$ are distinct and do not lie in $H$. 
\medskip

\noindent  If $\sum_{i=1}^m m_i \leq k$, then the dimension is  
$$N(k+l+2)+2k-4-\delta_{k,l} -2\sum_{i=1}^m m_i + m + \alpha, \ 
\mbox{where} $$

(i) $\alpha = 0$ if $C$ and $D$ both coincide with the directrix or are
contained in $H$;

(ii) $\alpha = n -\sum_{i=1}^{m+n} n_i$ if one of $C$ or $D$
coincides with the directrix and the other lies in $H$;

(iii) $\alpha = 2d -k-l+1$ if $C$ and $D$ coincide but are distinct
from the directrix and do not lie in $H$;

(iv) $\alpha = 2d -k -l + 1+ n  - \sum_{i=1}^{n+m} n_i  $ if $D$
coincides with the directrix or lies in $H$ and $C$ is a section
distinct from them;

(v) $\alpha =2d+2e-2k-2l+2+n -\sum_{i=1}^{n+m} n_i $ if $C$ and $D$
do not coincide; are distinct from the directrix and do not lie in
$H$.
\end{proposition}

\noindent {\bf Proof:}  The map from $F_{l-k}$ is given by $N+1$ sections
$s_0,\cdots,s_N$ in the class $\O_{F{l-k}}(e+lf)$. We assume $s_0$
corresponds to the hyperplane $H$. In case the sections $C$ and $D$ do
not coincide, we choose the $n$ lines along which $C$ and $D$ meet.
Since there is a one parameter family of fibers this gives us a choice
of $n$ dimensions. Finally, we mark two sections $C$ and $D$ on the
surface which have the required incidence with each other along the
chosen fibers. We need to know the dimension of pairs of sections with
specified incidence along specified fibers.

\begin{lemma}\label{independence}
  On $F_{l-k}$ the projective dimension of pairs of distinct
  irreducible sections $(s_1,s_2)$ in the classes $e+m_1f, e+m_2f$,
  $m_1 \leq m_2$, respectively, having contact of order $n_i$, $\sum_i
  n_i \leq m_1 + m_2 +k-l$, with each other along specified fibers is

\begin{displaymath}
2m_1+2m_2+2k-2l+2 - \sum_{i} n_i. 
\end{displaymath}
unless $m_1 = 0$ and $k \not= l$. In the latter case the dimension is 

\begin{displaymath}
2m_2 + k - l + 1 - \sum_i n_i.
\end{displaymath} 
\end{lemma}   

\noindent {\bf Proof:}  Set $u=m_1+m_2+k-l$ and $t=2m_2
+k-l+1$. We can pick the section in the class $e + m_1 f$ arbitrarily.
By Lemma \ref{dimsectionclasses} the dimension of these sections is
$2m_1 +k-l+1$, unless $m_1 = 0$ and $k \not= 0$. In the latter case
the dimension is $0$. We can embed the surface $F_{l-k}$ with the
section class $e+m_2f$.  In this embedding the first chosen section
class is a rational normal curve of degree $u$. To choose a curve in
the second section class with certain tangency conditions at specified
points is equivalent to choose a hyperplane in $\P^t$ with the
required tangency conditions to the rational normal curve along the
points of its intersection with the specified fibers. Since the curves
are distinct and irreducible the hyperplane will not contain the first
curve.  Hence, its intersection with the span of the first curve will
also be a hyperplane.

A rational normal curve is embedded by a complete linear system on
$\P^1$. Any two positive divisors of degree $d$ on $\P^1$ are linearly
equivalent. A hyperplane having
the specified order of contact is given by $n_1 p_1 + \cdots + n_r p_r
+ D$ where $D$ is any sum of points $m_1 + k_1 - \sum_i n_i$. Such
hyperplanes form a $u - \sum_i n_i$ dimensional linear subspace of
$(\P^{u})^*$. For each such linear space the hyperplanes of $\P^{t}$
that contain it is an irreducible linear space of dimension $m_2 -
m_1+1$ in $(\P^{t})^*$. The lemma follows.  $\Box$ \medskip

Using Lemma \ref{independence} we can complete the proof of
Proposition \ref{nastydim}. Note that if $\sum_{i=1}^m m_i > k$, then
the directrix of the surface has to be contained in $H$ and we must
have $\sum_{i=1}^m m_i = l$. In this case the dimension for the choice
of $s_0$ is $m+1$. By Lemma \ref{dimsectionclasses} the dimension for
the choice of each of the other $N$ section classes is $k+l+2$.
Finally, we have to choose the sections $C$ and $D$. 

\noindent $\bullet$ If they both coincide with the directrix,
there is nothing to choose.

\noindent $\bullet$ If they coincide, but are in a different
section class, by Lemma \ref{dimsectionclasses} we should add
$2d-k-l+1$.

\noindent $\bullet$ If $C$ coincides with the directrix, then we
have nothing further to choose. If $D$ does not coincide with the
directrix, then by Lemma \ref{independence} the choice for its
dimension is $2e-k-l+1- \sum_{i=1}^{n+m} n_i$.

\noindent $\bullet$ If $C$ and $D$ are distinct section classes
different from the directrix, Lemma \ref{independence} provides us the
dimension.
\smallskip

If $\sum_{i=1}^m m_i \leq k$, then the hyperplane in $H$ contains a
section class different from the directrix. The class in $H$ residual
to the lines must be $e+(l- \sum_{i=1}^m m_i) f$. Hence, by Lemma
\ref{dimsectionclasses} the dimension of $s_0$ is $k+ l - 2
\sum_{i=1}^m m_i + m +1$. The rest of the calculation is analogous to
the previous case. We must also projectivize and subtract the
dimension of the automorphism group of $S_{k,l}$. The proposition
follows. $\Box$ \medskip

\noindent {\bf Gluing scrolls.} We now prove that gluing
scrolls along fibers to form a tree imposes the expected number of
conditions. Since the fibers are $\P^1$ there is a three
dimensional family of isomorphisms which glue the two fibers. We will
often have to match section classes in different surfaces along the
fibers. This will impose additional conditions on the choice of
isomorphism.
\smallskip

Let $X \subset \P^N$ be a projective variety. Let $F_X(s)$ be an
irreducible, reduced subscheme of its Fano scheme of $s$ dimensional
linear spaces. Let $(p_i)_{i=1}^q$ denote $q \leq s+1$ points in
general linear position. Let $H(X)$ denote the $\P GL(N+1)$ orbit of
$[X]$ in the Hilbert scheme. For any $[Z] \in H(X)$, let $F_Z (s)$ be
the translation of $F_X(s)$ in the Grassmannian by the element that
takes $X$ to $Z$.

\begin{lemma}\label{superglue}
Let $X_1,X_2 \subset \P^N$ be two projective varieties. In the incidence
correspondence 
\begin{equation*}
\begin{split}
  & I := \{ (Z_1,\lambda_1, (p_j^1)_{j=1}^q, Z_2, \lambda_2,
  (p_j^2)_{j=1}^q): Z_i \in H(X_i), \lambda_i \in F_{Z_i}(s), \\ &
  p_j^1 \in \lambda_1, p_j^2 \in \lambda_2 \} \subset H(X_1) \times
  \G(s,N) \times (\P^N)^q \times H(X_2) \times \G(s,N) \times (\P^N)^q
\end{split}
\end{equation*}
the locus $J$ defined by $\lambda_1 = \lambda_2$ and $p_j^1 = p_j^2$
for all $1 \leq j \leq q$ has codimension $(s+1)(N-s) + sq$.

\end{lemma}

\noindent {\bf Proof:} $J$ is contained in the locus where $\lambda_1
= \lambda_2$. Restricted to this locus $J$ can be seen as the
projection that forgets the points on $\lambda_2$. Since the choice of $q$
points in $\P^s$ has dimension $sq$, the fibers of this projection are
all irreducible of dimension $sq$. Hence, to conclude the lemma it
suffices to settle the case $q=0$.  Consider the
projection from $I$ to the product of the Grassmannians $\G(s,N)
\times \G(s,N)$. By assumption all fibers are projectively equivalent
under the diagonal action of $\P GL(N+1)$, so they all have the same
dimension. We conclude that the codimension of $J$, which is the pull
back of the diagonal, is equal to the codimension of the diagonal in
$\G(s,N) \times \G(s,N)$.  $\Box$
\medskip

\noindent {\bf Forming trees.} Let $T$ be a connected tree with $v$
vertices numbered from $1$ to $v$.  Suppose $e_1$ of its edges are labeled
$1$, $e_2$ of its edges are labeled $2$ and the rest are labeled $0$. Let
$\nu_i$ be the valence of the $i$-th vertex.

Let $V$ be a $v$-tuple of scrolls where on the $i$-th one there are
$\nu_i$ marked lines with two distinct points on each. Let $W$ be the
subscheme of $V$ corresponding to the $v$-tuples that form the tree $T$.
That is if vertex $i$ is adjacent to vertex $j$, then a marked line of the
scroll in the $i$-th position coincides with a marked line of the scroll
in the $j$-th position and if the edge between the vertices is labeled
by $1$ or $2$, one or both of the points, respectively, coincide.
 
\begin{proposition}\label{glue}
The codimension of $W$ in $V$ is
\begin{displaymath}
(2N-2)(v-1) + e_1 + 2e_2.
\end{displaymath}
\end{proposition}

\noindent {\bf Proof:} When there are only two surfaces, the
proposition is a special case of Lemma \ref{superglue}. To prove the
proposition in general induct on the number of vertices by removing a
root and using Lemma \ref{superglue} again. Note that the same
argument applies when the lines joining a surface to two adjacent
surfaces coincide. $\Box$
\medskip

\noindent {\bf Main Enumerative Theorems:} We now state and prove the
main enumerative theorems under the assumption that the surfaces,
their limit hyperplane sections and the marked points on the curves
remain non-degenerate.
\smallskip

For an open subset of $\overline{\mathcal{M}}$ the image of the maps
have the same Hilbert polynomial. This induces a rational morphism to
the Hilbert scheme. The subspaces $\mathcal{X}$ and $\mathcal{W}$
whose general points correspond to maps with non-degenerate image also
admit rational maps to the same Hilbert scheme.

\begin{definition}
  A divisor of $\overline{\mathcal{M}}$ whose general point
  corresponds to a map with non-degenerate image is called {\bf
    enumeratively relevant} if its image in the Hilbert scheme has
  codimension one in the image of $\overline{\mathcal{M}}$.
\end{definition}

To determine the characteristic numbers of scrolls we count the number
of reduced points in the locus of the Hilbert scheme corresponding to
the type of scroll we are interested in. General one parameter
families of conditions lead to one parameter families of solutions. We
analyze the boundary divisors that these one parameter families meet.
Viewed from this perpective only enumeratively relevant divisors of
$\overline{\mathcal{M}}$ contribute to the enumerative calculations.

Contracted components can add moduli to maps from trees of Hirzebruch
surfaces without changing the number of moduli of the image surfaces.
In order to eliminate these extra moduli we define the modified tree.
\smallskip

\noindent {\bf The modified tree.} Let $\pi$ be a map from a tree
$\tilde{T}$ of Hirzebruch surfaces. We will refer to vertices
corresponding to surfaces that are contracted by $\pi$ as {\it
  contracted vertices} and to subtrees consisting entirely of
contracted vertices as {\it contracted subtrees}.  The {\it modified
  tree} $T$ is defined to be the same as $\tilde{T}$ if the map $\pi$
does not contract any surfaces. If $\pi$ contracts some surfaces, we
modify $\tilde{T}$ as follows:

\noindent $\bullet$ If a maximal, connected,
contracted subtree abuts only one non-contracted vertex, we simply
remove it.  

\noindent $\bullet$ If a maximal, connected, contracted subtree abuts exactly
two non-contracted vertices, we remove the subtree and join the two
vertices by an edge. 

\noindent $\bullet$ If a maximal, connected contracted subtree abuts
three or more non-contracted vertices, then we remove the contracted
tree and join the adjacent vertices to a node. We mark the node by
$c$, the number of non-contracted vertices adjacent to the contracted
tree. The reader can think of the node as a fiber common to more than
two surfaces in the image of $\pi$. We will call such lines
junctions (Figure \ref{tree}).

\begin{figure}[htbp]
\begin{center}
\epsfig{figure=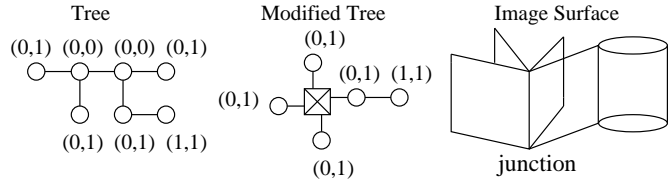}
\end{center}
\caption{The modified tree.}
\label{tree}
\end{figure}

We call a subtree of a modified tree connected, if every vertex has an
edge to some other vertex in the subtree or connects to a node that
another vertex in the subtree is connected to. The image of a
connected subtree is connected in codimension 1.  \smallskip

\noindent Although Theorems \ref{main1} and \ref{main2} initially look
long and complicated, in fact they state that only very few types of
behavior occur when we carry out the degenerations. Making all the
behavior precise, unfortunately, requires many cases. Colloquially,
the theorems assert that when we specialize a linear space $\Delta^I$
to a hyperplane $H$, some of the balanced scrolls satisfying the
incidences can remain outside $H$ (then the hyperplane section of the
scroll in $H$ meets $\Delta^I$), some can lie in $H$ (the scroll meets
the intersection of the linear spaces with $H$) and some can break
into a union of two balanced scrolls one of which is contained in
$H$---provided the surface and the limit of the hyperplane sections in
$H$ remain non-degenerate. In addition, the theorems specify that the
scrolls can become reducible only after there are enough linear
conditions to determine the limit of the hyperplane sections in $H$
and that $\Delta^I$ needs to meet the component contained in $H$.

\begin{theorem}\label{main1}
  Every enumeratively relevant component of the divisor
  $$\mathcal{D}_H \subset \overline{\mathcal{MS}}_H(N;
  k,l,C(k+l),D,Y,I,J_H)$$
  whose general member corresponds to a map
  where the set theoretic images of the surface and $C$ have maximal
  span is one of the following

\begin{equation*}
\begin{split}
  &1.\  \overline{\mathcal{MS}}_H (N; k,l,C(k+l),D, Y, I-1, J_H + 1:
  \Gamma^{J_H + 1} = \Delta^I) \\
  &2. \   \overline{\mathcal{MS}}_H(N-1;k,l,C(k+l),D,
  \tilde{Y},\tilde{I}, J_H : \Lambda^{\tilde{i}} = \Lambda^i \cap H,
  (\Delta^{\tilde{i}} = \Delta^i \cap H)_{i=1}^{I-1}, \\ & \ \ \ 
  \Delta^{\tilde{I}}=\Delta^{I}, J_H \geq
  k+l+1)   \\
  &3. \  \mathcal{X} (\P^N; (k_0, l_0; C(k_0 + l_0 +1),D(k_0): I \in
  \Delta(0), J_H (0) \geq k_0 + l_0 +2), \\ & \ \ \ (k_1= k-k_0, l_1=
  l-l_0, C(k_1+l_1-1),
  D(k_1))) \\
  &4. \  \mathcal{W} (\P^N; (k_i,l_i;C(d_i),D(k_i))_{i=0}^M : k_0 = 0,
  d_0 = k_0 + l_0 - M, \   \mbox{for} \  i > 0  \\ & \ \ \ d_i = k_i + l_i
  +1, \ I \in \Delta(i), \ J_H (i) \geq k_i + l_i +2  ).
\end{split}
\end{equation*} 
In addition, if $k=l$ there can be
\begin{equation*}
\begin{split}
 &5. \  \mathcal{X} (\P^N; (k_0, k_0+1; C(2k_0 + 2),D(k_0+1): I \in
  \Delta(0), J_H (0) \geq 2k_0 + 3), \ \ \ \  \\ & \ \ \ (k-k_0-1, 
  k-k_0, C(2k-2k_0-1),
  D(k-k_0))). \\
&6. \ \mathcal{X} (\P^N; (k_0, k_0+1; C(2k_0 + 2),D(k_0): I \in
  \Delta(0), J_H (0) \geq 2k_0 + 3), \\ & \ \ \ (k-k_0-1, 
  k-k_0, C(2k-2k_0-1),
  D(k-k_0+1))). 
\end{split}
\end{equation*}
\end{theorem}

\begin{theorem}\label{main2}
  Let $k+l \leq d \leq k+l+1$.  Every enumeratively relevant component
  of $$\mathcal{D}_{\Pi} \subset \overline{\mathcal{MS}}_{\Pi}
  (\P^N;k,l,C(d),D,Y,I,I',J_{\Pi})$$
  where the curve and the limit
  divisor cut out on $C$ by $\Pi$ remains non-degenerate is of the
  form
\begin{equation*}
\begin{split}
  &1. \  \overline{\mathcal{MS}}_{\Pi} (\P^N;k,l,C(d),D,Y,I,I'-
  1,J_{\Pi}+
  1: \Omega^{J_{\Pi} + 1} = \Sigma^{I'}) \\
  &2. \  \overline{\mathcal{MS}}_{\Pi} (\P^N;k,l,C(d),D,Y,I,I'-
  1,J_{\Pi}: q_{I'}' = p_j,
  \dim \Omega^{j} = \dim \Sigma^{I'} = N-2) \\
  &3. \  \overline{\mathcal{MS}}_{\Pi} (\P^N;k,l,C(d),D,Y,I,\tilde{I}'-
  1,J_{\Pi}+ 1: \pi(C) \subset \Pi, \\ & \ \ \ J_{\Pi} =d \leq k+l, \ 
  \Sigma^{\tilde{I}'} =
  \Sigma^{I} \cap \Pi)  \\
  &4. \  \overline{\mathcal{MS}}_{\Pi}
  (\P^{N-1};k,l,C(d),D,\tilde{Y},\tilde{I},\tilde{I}'- 1,J_{\Pi}+ 1:
  \pi(C) \subset \Pi, J_{\Pi} =d, \\ &\ \ \ \Delta^{\tilde{i}} =
  \Delta^i \cap \Pi,  \Sigma^{\tilde{i}'} = \Sigma^{i'} \cap \Pi)  \\
  &5. \  \overline{\mathcal{MS}}_{\Pi} (\P^{N};k,l,C(d),D,Y,I,I'-
  1,J_{\Pi}+ 1: C= \tilde{C} \cup F, \ 
  \pi(\tilde{C}) \subset \Pi, \\ & \ \ \  J_{\Pi} =d, q_{I'}' \in \tilde{C}) \\
  &6. \ \overline{\mathcal{MS}}_{\Pi} (\P^{N};k,l,C(d),D,Y,I,I'-
  1,J_{\Pi}+ 1: C= \tilde{C} \cup F, \ \pi(F) \subset \Pi, \\ & \ \ \ 
  J_{\Pi} \geq 2, p_{j_1}, p_{j_2},
  q_{I'}' \in F) \\
  &7. \ \overline{\mathcal{MS}}_{\Pi} (\P^{N};k,l,C(l),D,Y,I,I'-
  1,J_{\Pi}+ 1: C= D \cup F_1 \cdots \cup F_{l-k}, \\ &\ \ \ q_{I'}'
  \in \pi(D
  \cup \F_1 \cdots \cup F_r) \subset \Pi, \  J_{\Pi} \geq k+r+1) \\
  &8. \  \overline{\mathcal{MS}}_{\Pi} (\P^{N};k,l,C(l),D,Y,I,I'-
  1,J_{\Pi}+ 1:k>0,\ C= D \cup F_1 \cdots \cup F_{l-k}, \\ &\ \ \ 
  q_{I'}' \in \pi(\F_1 \cup \cdots \cup F_r) \subset \Pi, \ J_{\Pi}
  \geq 2r)
  \\
  &9. \ \mathcal{X} (\P^N; (k_0, l_0; C(k_0 + l_0 +1),D(k_0): I \in
  \Delta(0), J_H (0) \geq k_0 + l_0 +2), \\ & \ \ \ (k_1= k-k_0, l_1=
  l-l_0, C(d-k_0-l_0-1),
  D(k_1))). \\
  \end{split}
\end{equation*}
In addition if $k=l$ we can have 
\begin{equation*}
\begin{split}
 &10. \  \mathcal{X} (\P^N; (k_0, k_0+1; C(2k_0 + 2),D(k_0+1): I' \in
  \Delta(0), J_H (0) \geq 2k_0 + 3), \ \ \ \ \ \\ & \ \ \ (k-k_0-1, 
  k-k_0, C(d-2k_0-2),
  D(k-k_0))). \\
&11. \ \mathcal{X} (\P^N; (k_0, k_0+1; C(2k_0 + 2),D(k_0): I' \in
  \Delta(0), J_H (0) \geq 2k_0 + 3), \\ & \ \ \ (k-k_0-1, 
  k-k_0, C(d-2k_0-2),
  D(k-k_0+1))).
\end{split} 
\end{equation*}
\end{theorem}

\noindent We will prove both theorems by dimension counts. The
arguments become fairly intricate because we need to account for
competing phenomena.  In a tree of scrolls requiring three surfaces to
meet along the same line costs dimension. On a scroll requiring a
curve to contain a line also costs dimension. However, if a curve
contains the common line of three surfaces, then that voids the
matching conditions between pieces of a curve on different components
of a surface. The task at hand is to prove that the gain is always
less than the cost.  \smallskip

\noindent {\bf Proof:} 
By repeatedly using Propositions \ref{usekl} and \ref{linekl} and the
remark following them, it suffices to prove Theorem \ref{main1} when
$I=1, Y=0$ and Theorem \ref{main2} when $I=Y=0, I'=1$. The general
case then follows by adding marked points and marked fibers and
requiring them to lie in the intersection of general hyperplanes.
\smallskip

\noindent {\bf Proof of Theorem \ref{main1}:} By Proposition
\ref{nastydim} the dimension of $\overline{\mathcal{M}}_H$ is 
\begin{displaymath}
N(k+l+2) + 2k -4 + (J_H +2).
\end{displaymath}
The last term $J_H +2$ corresponds to the choice of $J_H$ points on
$C$ and a point on the surface. We can ignore these terms until the
end of the calculation, where they will help select the enumeratively
relevant divisors.

Let $U$ be a component of $\mathcal{D}_H$. Since we are interested
only in the enumeratively relevant components, we can work with the
modified tree $T$.  \smallskip

\noindent {\bf I. No component lies in $H$.} In this case, the hyperplane
section meets $\Delta^I$, i.e. $q_I \in C$. Hence the last term
decreases by 1 to $J_H +1$.  Since a reducible surface or a more
unbalanced surface that is the limit of scrolls $S_{k,l}$ has
dimension at most $N(k+l+2) + 2k -5$, the only divisor where the map
$\pi$ does not map a component of the surface to $H$ is given by 1.
\smallskip

\noindent {\bf II. The entire surface lies in $H$.} Since any limit of
scrolls that lies in $H$ is the limit of scrolls that lie in $H$. In
this case the dimension of the surface with a choice of hyperplane
section and directrix is at most $N(k+l+2) + 2k -5$ with equality if
and only if the surface is an $S_{k,l}$. Since the additional term we
add is still $J_H + 2$ we conclude that the only divisors has the
form 2. Furthermore, suppose $J_H < k+l+1$. Then there is a positive
dimensional space of hyperplane sections that pass through the chosen
$J_H$ points. Such a divisor is not enumeratively relevant. 
\smallskip

\noindent {\bf III.} From now on we can assume that $\pi$ maps at
least one component of the surface into $H$ and keeps at least one
component outside $H$.  \smallskip

\noindent {\bf Using the non-degeneracy assumption to simplify $T$.}  Since
we are assuming that the image surface spans $\P^{k+l+1}$, each of the
non-contracted components  map to rational normal scrolls.
Moreover, the span of any subsurface connected in codimension 1 of
degree $d$ must be $\P^{d+1}$. Since the limit curve $C$
is non-degenerate each of its non-contracted components maps to a
rational normal curve. We conclude that the restriction of $C$ to any
subsurface of degree $d$ connected in codimension 1 has degree at most
$d+1$ since such a surface can span at most $\P^{d+1}$.
\smallskip
 
{\it The surfaces outside form a connected tree and have multiplicity
  1 with $H$ along their lines contained in $H$.} Let $T_i$ be a
maximal, connected subtree of the modified tree $T$ where all the
vertices correspond to surfaces mapped into $H$. At most one surface
outside $H$ can be connected to $T_i$. Otherwise, the restriction of
$C$ would have degree at least 2 more than the total degree of the
surfaces in $T_i$ contradicting the previous observation.  As a
corollary we conclude that the surfaces outside $H$ form a single
connected tree.  Again by the non-degeneracy of $C$ the surfaces
outside have multiplicity 1 along any line they meet $H$ including the
lines where they are connected to surfaces inside. Moreover, a line
joining two surfaces outside $H$ is not contained in $H$.  \smallskip

\noindent {\bf Reminder:}  We remind the reader that a line of the
image surface common to more than two components is called a {\it
  junction}. A line where exactly two components meet is an {\it
  ordinary common line}.

\smallskip

\noindent {\bf Bounding the dimension of maps.} We now calculate the
dimensions of various loci of maps.

{\bf Notation:} We assume that $\pi$ is a  map whose image forms a
contracted tree $T$ of scrolls with $v$ vertices and $\sigma$ nodes
marked $c_1, \cdots, c_{\sigma}$. We assume that $T$ has 

 $\bullet$  $t$ maximal connected subtrees $T_1, \cdots, T_t$
with $v_1, \cdots, v_{t}$ vertices whose image under $\pi$ lies in $H$
and

 $\bullet$ one connected subtree $T_0$ with $v_0$ vertices
all of which lie outside $H$. Let the total number of vertices be $v$.
\smallskip

\noindent We assume that the $i$-th vertex is marked by two integers
$(k_i,l_i)$ to denote that the surface corresponding to that vertex is
$S_{k_i,l_i}$. We must, of course, have $\sum k_i + l_i = k+l$. In
addition, each vertex has the information of two limit curves $C,D$
associated to it. We assume that 

 $\bullet$ Both $C$ and $D$ contain $\sigma_1$ junctions
$F_1, \cdots, F_{\sigma_1}$ and $\mu_1$ ordinary common lines $G_1,
\cdots , G_{\mu_1}$ where $D$ contains them with multiplicity $f_i$
and $g_i$, respectively,

 $\bullet$ $D$ contains $\mu_2$ additional ordinary common lines and
$\sigma_2$ junctions with multiplicities $f_i$ and $g_i$,
respectively, and

 $\bullet$ $C$ contains $\mu_3$ additional ordinary common lines
and $\sigma_3$ junctions.

$\bullet$ On the $j$-th surface apart from the ordinary common
lines and nodes, $C$ and $D$ contain $\zeta_j$ common fibers, $C$
contains $\xi_j$ additional fibers and $D$ contains $\psi_j$
additional fibers. We allow $D$ to contain the common fibers with
multiplicity $\rho_{ji}$ and the additional fibers with multiplicity
$\phi_{j,i}$.

$\bullet$ Finally, on surfaces where the section parts of
$C$ and $D$ are distinct, the sections have contact of order $m_{i,j}$
with each other along fibers common to at least two surfaces and
contact of order $n_{i,j}$ along other fibers common to both $C$ and
$D$.  \smallskip

\noindent {\bf Caution:}  Recall our convention that {\it  a perfectly
  balanced scroll does not have any directrices.}  \smallskip

\noindent {\bf Summation convention.} $\sum_{EE}, \sum_{CC},
\sum_{CE}, \sum_{ED}, \sum_{CD}$ will denote summation over indices
for which both $C$ and $D$ contain the directrix, both $C$ and $D$
contain the same section class, $D$ contains the directrix, $C$
contains the directrix, or they both contain different section
classes, respectively. We will also decorate the summation notation
with $I$ or $O$ to denote summing over only the surfaces inside or the
surfaces outside.
\smallskip

\noindent {\bf The dimension of the building blocks.} We have to
arrange our surfaces outside so that their intersection with $H$ is
governed by the data of $C$. In addition, we have to choose $D$ in all
the surfaces and $C$ in the surfaces inside. Using Proposition
\ref{nastydim} we see that the dimension of the data of the surfaces
$S_{k_i,l_i}$ together with the curves $C$ and $D$ on them prior to
the gluing conditions is

\begin{equation*}
\begin{split}
  &\sum_{j=1}^v (N(k_j + l_j +2) + k_j - l_j -6 - \delta_{k_j,l_j}) +
  \sum_{EE|O} (l_j +1)+ \sum_{CC|O} (d_j +2) \\ &+ \sum_{ED|O} (2e_j -k_j
  +2) + \sum_{CE|O} (d_j +2) + \sum_{CD|O} (d_j +2e_j -k_j-l_j + 3) +
  \sum_O \psi_j \\ &+ \sum_{CC|I} (2d_j -k_j -l_j + 1) + \sum_{ED|I} (2e_j
  -k_j-l_j +1) + \sum_{CE|I} (2d_j -k_j-l_j +1) \\ &+ \sum_{CD|I}2 (d_j
  +e_j - k_j-l_j +1) +\sum_{I} (\zeta_j + \xi_j + \psi_j) -
  \sum_{ED,CE \atop CD}\sum_i (m_{ji}+n_{ji}) 
\end{split}
\end{equation*}
{\bf Gluing the pieces.} So far we have the dimension of $v$ surfaces
$S_{k_j,l_j}$, all but $v_0$ of them lying in $H$ and the choices of
two section classes on each individual surface. We have to glue the
surfaces along fiber lines and we have to make sure that both $C$ and
$D$ are connected.  We have to match the points where the section
parts of $C$ and $D$ meet the common fibers of two surfaces, except
when the common fiber is contained in $C$ or $D$. In the latter case
connectedness of the curve is automatic. When a fiber joining the
surfaces is not contained in the curves, we normally have two sections
to match, except when the sections already meet the fiber at the same
point, i.e. $m>0$. We now introduce some more notation to record all
these different possibilities.  \smallskip

\noindent {\bf More Notation.} Note that $T$ has
$v-1-\sum_{i=1}^{\sigma}(c_i-1)$ ordinary common lines. Suppose

$\bullet$ in  $\alpha_1$ of them surfaces with distinct sections meet and
$C$ and $D$ both contain the common fiber; in $\alpha_2$ of them
distinct surfaces with distinct sections meet and only $D$ contains
the fiber; and in $\alpha_3$ of them surfaces with distinct sections
meet and $C$ contains the fiber;

$\bullet$ in $\beta_1$ of them surfaces with common sections meet and
both $C$ and $D$ contain the fiber; in $\beta_2$ of them only $D$ and
$\beta_3$ of them only $C$ contains the fiber.

Note that there must be $\mu_1 - \alpha_1 - \beta_1$ ordinary common
lines where a surface with distinct sections meet a surface with
common sections and both $C$ and $D$ contain the fiber. Similarly, for
the other two cases. 

Other than these fibers suppose that there are $\kappa_1$ ordinary
fibers where surfaces with 
distinct sections meet and have $m_{ji} > 0$, $\kappa_2$ ordinary
fibers where surfaces with distinct sections meet and
$m_{ji}=0$. $\kappa_3$ fibers where surfaces with distinct sections
meet those with common sections and $\kappa_4$ fibers where surfaces
with common sections meet. We remark that 
\begin{displaymath}
v-1-\sum_{i=1}^{\sigma}(c_i-1) = \sum_{i=1}^3 \mu_i + \sum_{i=1}^4 \kappa_i 
\end{displaymath}

Suppose at each of the first $\sigma_1 + \sigma_2 + \sigma_3$ junctions where
either $C$ or $D$ or both contain the fiber $\tilde{c_i}$ surfaces
with common sections meet $c_i'$ surfaces with distinct
sections. Suppose in the next $\sigma_4$ junctions $\tilde{c_i}$ surfaces
with common sections meet $c_i'$ surfaces with distinct sections. Of
course, we must have $c_i = \tilde{c_i} + c_i'$. In the next
$\sigma_5$ junctions suppose that surfaces with distinct sections, but
with $m>0$ meet. Finally, in the last $\sigma_6$ junctions surfaces with
distinct sections and $m =0$ meet. Note that $\sum_{i=1}^6 \sigma_i = \sigma$.
\smallskip

\noindent {\bf Choosing the lines to glue along.} On each surface we
have to choose the lines along which we will glue it to other
surfaces. For gluing two surfaces outside or two surfaces inside this
gives us a choice of $2$ for each connection except at the junctions
where this over counts by $c_i -2$. For gluing a surface outside and a
surface inside we have only a choice of 1 dimension because the
surface outside has only finitely many lines in $H$.
\smallskip

\noindent {\bf Combining the gluing conditions.} When we glue the
sections if in one surface the sections coincide, then on the other
surface we must have $m>0$ or the two sections must coincide, hence we
only have 1 point to match. In case the sections do not match on
either of the two surfaces, we have two points to match, unless $m>0$
on one, hence both of them.  By Lemma \ref{superglue} and Proposition
\ref{glue} we conclude that we have to add the following gluing
conditions to our previous expression

\begin{equation*}
\begin{split}
&\sum_{i=1}^t (6-2N) (v_i -1) + (4 - 2N)(v_0 -1) + (5-2N) t +
\sum_{i=1}^{\sigma} (2-c_i) \\ &-\mu_2 - \mu_3 -\kappa_1 - 2 \kappa_2 -
\kappa_3 - \kappa_4 - \sum_{i=\sigma_1+1}^{\sigma -  \sigma_6}
(c_i -1) - \sum_{i=\sigma - \sigma_6 +1}^{\sigma} 2(c_i-1)
\end{split}
\end{equation*}   
{\bf Simplifying the formulae.} Now we simplify these expressions
using the facts that $D$ has degree $k$, $C$ has degree $k+l$ and that
the total surface has degree $k+l$. In the computations I found it
helpful to express these facts in terms of the following equations

\begin{equation*}                                                             
\begin{split}
  k &= \sum_{EE,CE} k_j + \sum_{CC,ED \atop CD} e_j + \sum_{i,j} (\rho_{ji}
  + \phi_{ji}) + \sum_{i=1}^{\sigma_1 + \sigma_2} f_i +
  \sum_{i=1}^{\mu_1 + \mu_2} g_i   
 \\ k+l &= \sum_{ED,EE} k_j + \sum_{CE,CC \atop CD} d_j + \sum_j (\zeta_j +
  \xi_j) + \mu_1 + \mu_3 + \sigma_1 + \sigma_3 
  \\ k+l &= t + \sum_I (k_j +l_j) + \sum_O (\zeta_j + \xi_j) +
  \sum_{EE,ED|O} k_j + \sum_{CE,CC|O \atop CD} d_j 
\end{split}                                                                  
\end{equation*}   
Finally, we have to use the fact that the curves cannot have total
intersection multiplicty larger than $k$ at isolated points on the
smooth locus of the surface. If we assume that the total intersection
multiplicity is $k-w$ for some nonnegative $w$, we obtain the
relation:

\begin{equation*}
\begin{split}
  k-w &= \sum_{ED} (e_j+ \xi_j + \phi_j -l_j) + \sum_{CE} (d_j+
  \xi_j + \phi_j -l_j) \\ &+ \sum_{CD} (e_j+ \xi_j + \phi_j + d_j -k_j
  -l_j) - \sum_j \left( \sum_i m_{ji}+\sum_i n_{ji} \right)
\end{split}
\end{equation*}
{\bf Final formula.} It is straightforward, but messy algebra to
conclude that the locus of maps with the tree we described has
dimension at most

\begin{equation*}
\begin{split}
  &N(k+l+2) + 2k - 3 - v_0 - w + \sum_{CC,CD} (k_j - l_j)- \sum_j
  \delta_{k_j,l_j} -\sum_{EE} 1 + \sum_{CD} 1 \\ & - \sum_{ED, CE
    \atop CD} (\xi_j + \sum_i \phi_{ji}) + \sum_j \psi_j - \sum_j
  \sum_i (\rho_{ji} + \phi_{ji}) - \sum_{i=1}^{\sigma_1 + \sigma_2}
  f_i -
  \sum_{i=1}^{\mu_1+\mu_2} g_i   \\
  & - \mu_3+ \sigma -\sigma_1 -\sigma_3 -
  \sum_{i=\sigma_1+1}^{\sigma - \sigma_6}(c_i -1) -
  \sum_{i=\sigma - \sigma_6 + 1}^{\sigma}
  2(c_i-1) - \kappa_2
\end{split}
\end{equation*}

\noindent {\bf Interpreting the formula.} We compare 
$N(k+l+2) + 2k -4$ with the formula to see when the formula is exactly
one less.

\noindent $\bullet$ Since  $k_j - l_j - \delta_{k_j,l_j} \leq -1$, the
surfaces of the form $CD$ contribute less than or equal to $0$ to the
sum with equality if and only if the surface is balanced.

\noindent $\bullet$  Note that $\psi_j - \sum_i \phi_{ji}
\leq 0$.

\noindent $\bullet$ For each junction between $1$ and $\sigma -
\sigma_6$ the contribution is $-1$ or less and for every other
junction the contribution is less than $-2$.

\noindent $\bullet$ All the  other terms after
$2k$, if they exist, are strictly negative. 

Since $v_0 \geq 1$,
$$N(k+l+2) + 2k -3 -v_0 \leq N(k+l+2) + 2k -4$$
we conclude that

\noindent $\bullet$ $\sigma_6 = 0$,

\noindent $\bullet$ There can be at most one surface of type $EE$ or
$CC$ since these surfaces contribute $-1$ or less to the sum,

\noindent $\bullet$ If the surface contains a subsurface of type $EE$ or
$CC$, $C$ and $D$ cannot contain any junctions, ordinary common
lines, or any fibers on surfaces of type $CE,ED,CD$.

Using the equations above we can reexpress $w$ as
\begin{equation*}
\begin{split}
  w &= \sum_{CC} (2d_j-k_j-l_j) + \sum_{EE} (k_j-l_j) + \mu_1 + \mu_3
  + \sigma_1 + \sigma_3 + \sum_j \zeta_j \\ &+\sum_{i=1}^{\sigma_1 +
    \sigma_2} f_i + \sum_{i=1}^{\mu_1+\mu_2} g_i +\sum_{CC,EE}
  (\xi_j + \phi_{ji}) + \sum (m_i +n_i)
\end{split}
\end{equation*}   
$\bullet$  A codimension $1$ locus satisfying our assumptions
does not contain any components of type $CC$. If there is a surface of
type $CC$, then $w$ is at least $1$ since $m \geq 1$ for the
components abutting the component of type $CC$ and $2d_j -k_j -l_j
\geq 0$.
\smallskip

\noindent {\bf Suppose there are no surfaces of type $EE$.} Then the
surfaces are of types $CE, ED$ or $CD$. Any ordinary common line or
junction contained in the curves contributes less than or equal to
$-2$. Hence, $C$ and $D$ cannot contain any junctions or ordinary
common lines. Each time we glue two of these surfaces either $\kappa_2
= 1$ or $m \geq 1$ for both of the surfaces. We conclude that at a
general point in a codimension one locus there must be exactly two
surfaces, one inside $S_{k_I,l_I}$ and one outside $S_{k_O,l_O}$, and
the curves $C$ and $D$ meet the common fiber at distinct points. Since
$C$ and $D$ cannot contain fibers (this would contribute less than or
equal to $-1$), their restriction to each component must be a section.
The degree of the component of $C$ in the surface outside is $k_O +
l_O - 1$ and the degree of the component in the surface inside is $k_I
+ l_I +1$ .

The surfaces have the types $(CD,CE)$, $(CE,ED)$ or $(CE,CE)$. 

\noindent $\bullet$ If the type is $(CE,CE)$, then $k= k_1 + k_2$.

\noindent $\bullet$ If a surface is of type $ED$, then that surface must
be outside $H$. Since the directrix of the surface has degree one less
than its degree the surface must be a plane. (Remember the quadric
surface does not have any directrices.)  The degree of $D$ is at least
$1$ and it is required to pass through a point. Hence the choice of $D$
contributes a one-parameter family to our dimension. Such a component
cannot be enumeratively relevant unless $k=l$. In the latter case a
perfectly balanced scroll breaks into a plane and another balanced
scroll. 

\noindent $\bullet$ If there is a surface of type $CD$, it must be
balanced. Either the scroll is perfectly balanced and $D$ is a
directrix, in which case $k = k_1 + k_2$ or the choice of $D$
contributes at least 1 to the dimension. It contributes exactly one
when the surface is balanced and not perfectly balanced and the degree
of $D$ is one larger than the degree of the directrix. Such a case can
be enumeratively relevant only if $k=l$. In that case we must have $k= k_1
+ k_2 +1$ and a perfectly balanced scroll breaks into a union of two
balanced scrolls.  \smallskip
\smallskip

\noindent {\bf Suppose that there is a component $S$ of type
  $EE$.} As we already observed in this case there are no junctions
and the curves $C,D$ do not contain any ordinary common lines or lines
on the surfaces of other types. Furthermore, we must have $w=0$.
Hence, the equation for $w$ simplifies to $$0= k_S - l_S + \xi_S +
\sum \phi_{ji} + \sum_i (m_{ji} + n_{ji}).$$
The degree of the
hyperplane section restricted to $S$ is $k_S + \xi_S$. Since any
surface emanating from $S$ can account for at most one degree by our
non-degeneracy assumption, we conclude that there must be at least
$l_S - \xi_S$ surfaces adjacent to $S$. They are surfaces of type
$CE,ED$ or $CD$, so $m \geq 1$ for each of them.  Since $w=0$, we
conclude that there are exactly $l_S - \xi_S$ surfaces, all adjacent
to $S$. The curves $C$ and $D$ meet each other simply along the lines
joining the other surfaces to $S$.  $S$ lies outside $H$ and all the
other surfaces are in $H$. Moreover, $k_S = 0$, so $S$ is a cone.

The surfaces inside have to be of type $CE$ or $CD$. If they are all
of type $CE$, then $k = \sum k_i$. If some surfaces are of type $CD$,
either all of them are perfectly balanced and $D$ is a directrix.
Hence, $k = \sum k_i$.  Or $k=l$ the surface is balanced, the degree
of $D$ is one more than the degree of the directrix. There can be at
most one such surface. In this case $S$ must be a plane and there are
only two components.

Finally, we have to add the choices for the $J_H$ points and the point
$q_I$. We must have $q_I$ in a surface inside, otherwise we lose 1
dimension. Finally, to be enumeratively relevant $C$ should not move
in a linear system. We conclude that $J_H (I) \geq k_I + l_I +2$ for
each of the surfaces inside. (In fact, $J_H(I) \geq k_I + l_I +3$ in
case the surface outside is a plane.) This completes the proof. $\Box$
\bigskip

\noindent  {\bf Proof of theorem \ref{main2}.} It suffices to prove
the case $I=Y=0, \  I'=1$. The dimension of
$\overline{\mathcal{MS}}_{\Pi}$ is
\begin{displaymath}
(N+1)(k+l+2) + 2d -2l -5.
\end{displaymath}
Let $U$ be a component of $\mathcal{D}_{\Pi}$. We have $J_{\Pi}
\leq d.$

\noindent {\bf 1.} Suppose the surface is irreducible. If the curve is
also irreducible, then $q_{I'}' \in C \cap \Pi$.  Here $q_{I'}'$ can
also coincide with one of the $p_j$. Of course, we mean that 
there is a map to $\F(0,0,1;N)$ with a contracted component and the
marked points move to the contracted component. We will continue using
the more suggestive notation $q_{I'}' = p_j$ since this is more
convenient for counting purposes.  If the curve does
not lie in $\Pi$, these are  divisors provided that the surface remains
an $S_{k,l}$. If the curve lies in $\Pi$, then either $d \leq k+l$ or the
surface also lies in $\Pi$. In either case this locus is
a divisor precisely when $J_{\Pi}=d$.

If the curve becomes reducible, then it must be the union of a section
of degree $d'$ together with $d-d'$ fibers. We denote this curve
as $\tilde{C} \cup F_1 \cup \cdots \cup F_{d-d'}$. If $\tilde{C}$
is contained in $\Pi$, then  the locus is a divisor when $d' = d-1$
and $J_{\Pi} = d$ and the remaining fiber is not contained in $\Pi$ or
$d=l$ and the directrix union $r$ fibers are contained in $\Pi$ and
$J_{\Pi} = k+r+1$.

If only fibers are contained in $\Pi$, then either there must be a
unique fiber and two points among the $J_{\Pi}$ must specialize to it
or $d=l$ and the curve must break to the union of the directrix with
lines $r$ of which are contained in $\Pi$.  In the latter case $2r$ of
the points in $J_{\Pi}$ must specialize pairwise to the lines..

\noindent {\bf 2.}  We can assume that both the curve and the surface
break. (We include the case when the surface breaks into planes unions
scrolls and the image of the curve does not `really' break.) We will
perform a calculation analogous to the previous case to find the
divisors. We work with the modified tree $T$.  \smallskip

\noindent {\bf Using the non-degeneracy assumption to simplify $T$.}
Since the hyperplane section of $C$ in $\Pi$ is non-degenerate, the
components of $C$ outside $\Pi$ meet $\Pi$ transeversely. Each tree of
curves inside can meet at most one curve outside. In particular, the
curves outside form a connected tree. The surfaces outside have
contact of order one with $\Pi$ along their lines contained in $\Pi$
except possibly when the surface is a cone and the section part of $C$
reduces to the vertex. However, in the latter case forcing the cone to
have higher order of contact with $\Pi$ strictly lowers the dimension,
so we can ignore this case.  \smallskip

\noindent {\bf Three or more curves.} Unfortunately, in  
addition to the cases we considered in the previous proof, now the
surfaces outside $\Pi$ can have sections that lie in $\Pi$ or remain
outside $\Pi$.  Fortunately, due to our non-degeneracy assumptions we
do not need to record the incidence data of $C$ with the hyperplane.
If we remove the non-degeneracy assumption, it is easy to construct
situations where the directrix has high order of contact with the
hyperplane and the curve has high order of contact with both the
hyperplane and the directrix at different points. It is a very hard
problem to determine the limits in this generality. This is the main
obstruction for carrying out our algorithm in general.

\smallskip

\noindent {\bf Notation:} We preserve the notation and conventions of
the previous proof. We decorate our sums for the surfaces outside $H$
with $* \subset H$ to signify summation over the incides where $*$ is in
$H$. We decorate the summation by $'$ if none of $E_j,C_j,D_j$ is
contained in $H$.  In that case, we let $x_j$ denote the number of
lines that the surface has in $H$. Finally, we decorate our lines with
$H$ and $\not\subset H$ in case the line belongs to the part of $C$ or
$D$ that lies in a surface outside $H$, and the line is in $H$ or not,
respectively.

\noindent {\bf Building blocks:}  By Proposition \ref{nastydim} the
dimension prior to gluing is

\begin{equation*}
\begin{split}
  &\sum_{j=1}^v (N(k_j + l_j +2) -6 + k_j -l_j -\delta_{k_j,l_j}) +
  \sum_{EE|O}' (k_j+l_j-x_j+2) \\ &+ \sum_{EE|O \atop E \subset H}
  (l_j + 1) + \sum_{CC, CE|O \atop C\subset H} (d_j + 2) + \sum_{CC,
    CE|O \atop E \subset H} (2d_j -k_j+2) + \sum_{ED|O \atop D \subset
    H} (e_j + 2) 
\end{split}
\end{equation*}
\begin{equation*}
\begin{split}
 & +\sum_{CC, CE|O}' (2d_j - x_j +3) + \sum_{ED|O
    \atop E\subset H} (2e_j -k_j +2) + \sum_{CD|O \atop C \subset
    H}(d_j + 2e_j -k_j-l_j +3) \\ &+ \sum_{CD|O \atop D \subset
    H}(e_j + 2d_j -k_j-l_j +3) + \sum_{CD|O \atop E \subset H}(2e_j +
  2 d_j -2k_j -l_j + 3) \\ &+ \sum_{CD|O}' (2d_j + 2e_j -k_j -l_j -x_j
  + 4)+ \sum_{ED|O}' (2e_j - x_j +3) \\& + \sum_{ED|I} (2d_j -k_j-l_j
  +1) + \sum_{CE|I} (2d_j -k_j-l_j +1) \\ &+ \sum_{CC|I} (2d_j -k_j
  -l_j + 1)+ \sum_{CD|I} 2 (d_j +e_j - k_j-l_j +1) + \sum_{I} (\zeta_j
  + \xi_j + \psi_j) \\ & + \sum_O (\zeta_{j,\not\subset H} + \xi_{j,
  \not\subset H} + \psi_{j, \not\subset H})
  - \sum_{ED,CE \atop CD}\sum_i (m_{ji}+n_{ji})
\end{split}
\end{equation*}

\noindent {\bf Final Formula.} The gluing conditions are identical to the
previous case. We retain the notation we used there. Finally we use
the facts that the degree of $C$ is $d$, the degree of $D$ is $k$ and
that their intersection multiplicities at isolated points along the
smooth points of the surface has to be $d-l-w$ for some positive $w$.
Combining everything we obtain

\begin{equation*}
\begin{split}
  & (N+1)(k+l+2) + 2d-2l-5 -v_0 - w - \sum_{E \subset H|O}k_j -
  \sum_{C \subset H |O} d_j - \sum_{ED,CD|O \atop D \subset H} e_j \\&
  - \sum_I (k_j+l_j) -t - \sum_j \delta_{k_j,l_j} + \sum_{CC,CD} (k_j
  -l_j) + \sum_I \psi_j + \sum_O \psi_{j, \not\subset H} - \sum_i
  \phi_{ji} \\& - \sum_{1}^{\sigma_1 + \sigma_2} f_i - \sum_1^{\mu_1
    +\mu_2} g_i - \sum_{j,i} \rho_{ji} - \mu_3 - \sum_{CD,CE \atop ED}
  (\xi_j + \sum_i \phi_{ji}) + \sigma -\sigma_1 -\sigma_3 \\ & -
  \sum_{\sigma_1+1}^{ \sigma - \sigma_6}(c_i -1) - \sum_{\sigma-
    \sigma_6 + 1}^{\sigma} 2(c_i-1)- \kappa_2 - \sum_{EE} 1 + \sum_O '
  (1 - x_j) + \sum_{CD} 1 \\& - \sum_O ( \psi_{j, \subset H} +
  \zeta_{j,\subset H} + \xi_{j,\subset H})
\end{split}
\end{equation*}

\noindent {\bf Interpretation.} This expression gives us the dimension
of the tree of scrolls before we choose the $J_{\Pi}$ points that are
the marked points of $C \cap \Pi$. Since we are assuming that the
points remain non-degenerate for each connected tree of curves in $H$
we can choose one more point than the total degree of the curves in
the tree. Once we choose the points we compare the above expression
with $(N+1)(k+l+2)+2d-2l-5$.

The curves in a tree of surfaces in $\Pi$ form a connected
tree. Consequently, there can be two types of connected trees of
curves in $\Pi$: Trees that contain a curve of a surface contained in
$\Pi$ and trees that do not. We refer to these trees as trees of the
first and second kind, respectively.

\noindent $\bullet$ The terms $\sum_I (k_j +l_j) + t$ and $d_j +
\zeta_{j,\subset H} + \xi_{j,\subset H}$ taken over the curves and
surfaces contained in a tree of the first kind add up to at least the
degree of the curve contained in it. Hence with the choice of points
the contribution of each such tree is at most 1 with equality if and
only if the degree of the curve in each tree of surfaces inside is the
maximum allowed. We can also assume that the trees on the surfaces in
$\Pi$ do not get connected because otherwise $t$ contributes
negatively. Note that junctions or ordinary common lines do not change
the conclusion and in fact continue to contribute negatively. We can
further assume that trees of the first kind do not have any curves or
fibers (except for junctions or ordinary common lines) contained in
surfaces outside $\Pi$. If such a curve exists, there must be an
additional surface outside contributing to $v_0$, hence bringing the
contribution of the tree to at most 0.

Surfaces outside have to connect these $t$ trees of curves. Let us
call such surfaces connecting surfaces. A connecting surface cannot
have a section part of $C$ in $\Pi$ since otherwise it would connect
the trees inside. If it has $D$ in $\Pi$, then $-e_j$ would contribute
an amount more than the number of trees it connects. If it has $E$ in
$\Pi$, the number of trees the surface connects would be at most one
more than $k_j$. Together with the contribution from $v_0$, this would
annul the contribution of all those trees. Finally, if the surface has
no special sections, then it must contain a line for every surface it
connects, hence together with $v_0$ the term $1-x_j$ annuls the
contribution of the trees it connects.

The argument for trees of the second kind is analogous but easier. We
conclude that the choice of points at most exactly cancel the
contribution of the terms discussed.

\noindent $\bullet$  If there is a surface outside which contains no
special sections and no lines in $\Pi$, the contributions of $v_0$ and
$(1-x_j)$ exactly cancel each other out.

\noindent $\bullet$ Since $k_j -l_j - \delta_{k_j,l_j} \leq -1$ a
component of type $CC$ contributes $-1$ or less and a component of
type $CD$ contribues $0$ or less to the sum with equality if and only
if they are balanced. A component of type $EE$ contributes $-1$. We
conclude that there can be at most one surface of type $EE$ or $CC$.

To continue our analysis we express  $w$ as in the previous
proof. By arguments analogous to the ones given there we have that 

\noindent $\bullet$ a codimension 1 locus does not contain any
components of type $CC$.

\noindent $\bullet$ if the surface does not contain a subsurface of type
$EE$, there can be at most 2 components. The curves $C$ and $D$ do not
contain any fibers and the degree of the curve in the surface inside
is $k_I + l_I +1$. The enumeratively relevant codimension one loci
have $k_1 + k_2 = k$, unless $k=l$ and $k_1 + k_2 = k-1$.

\noindent $\bullet$ a subsurface of type $EE$ is not contained in
$\Pi$. Suppose a surface contains a subsurface $S$ of type $EE$. There
cannot be any junctions or ordinary common lines in this case. All the
other surfaces are of type $CE,CD,$ or $ED$. By an argument similar to
the previous case, we see that there cannot be a component of type
$EE$ in case the curve has degree $k+l+1$. In case the degree of the
curve is $k+l$ there can be components of type $EE$ only in the case
described in the previous theorem. However, then all the curves inside
are connected and we conclude that there can be at most one other
component. One can continue the analysis to see what components would
appear if the degree of the curve were smaller. When one would like to
apply the algorithm in cases we will not explain in this paper, this
extension becomes useful. $\Box$

\section{Multiplicity Calculations}

In this section we carry out the multiplicity calculations needed for
enumerative computations involving balanced scrolls. The philosophy,
motivated by Vakil's work on curves, is that any multiplicity should
be reflected in the local structure of the limit surface.

Under the hypotheses of Theorems \ref{main1} and \ref{main2} when we
specialize a linear space to $H$ or $\Pi$, balanced scrolls incident
to the linear space break into at most two balanced scrolls. By
Proposition \ref{direc} the limit of the directrices is uniquely
determined by the surfaces. For multiplicity calculations it is more
convenient to reformulate the degeneration problem in the space of
scrolls where we only mark the hyperplane section.

Let $\overline{\mathcal{M}}_H ( \P^N; k,l;C, \{ \lambda_i
\}_{i=1}^Y,\{ q_i \}_{i=1}^I, \{ p_j \}_{j=1}^{J_H} )$, or
$\overline{\mathcal{M}}_H$ for short, be defined like the
corresponding space $\overline{\mathcal{MS}}_H$ in \S 6 except that
now do not mark the directrix. More explicitly, an open set in
$\overline{\mathcal{M}}_H$ corresponds to maps from a Hirzebruch
surface into $\P^N$ as a scroll $S_{k,l}$, where the marked curve $C$
maps to the hyperplane section in $H$ and the marked fibers and points
are required to lie in various linear spaces. We compactify the space
as in \S 3.  In a similar fashion define the space
$\overline{\mathcal{M}}_{\Pi}$ and the loci $\mathcal{X}_H'$ and
$\mathcal{X}_{\Pi}'$ corresponding to $\overline{\mathcal{MS}}_{\Pi}$
and $\mathcal{X}_H$ and $\mathcal{X}_{\Pi}$, where again the only
difference is that we do not mark the directrix. The spaces
$\overline{\mathcal{M}}_H$ and $\overline{\mathcal{M}}_{\Pi}$ have
natural Cartier divisors $D_H'$ and $D_{\Pi}'$ defined by requiring
the point $p_I$ and $q_{I'}'$ to lie in $H$ and $\Pi$, respectively.

Since the limits of directrices are determined uniquely by the
surfaces, Theorems \ref{main1} and \ref{main2} describe the
enumeratively relevant components of $D_H'$ and $D_{\Pi}'$ subject to
the non-degeneracy assumptions. We would like to compute the
multiplicity of the Cartier divisor along each of the Weil divisors
appearing in the list. \smallskip

\begin{lemma}\label{normal1}
  Let $S_{k,l}$ be a non-singular scroll in $\P^N$. Let $\nu :=
  \nu_{S_{k,l} / \P^N} $ denote its normal bundle in $\P^N$. Suppose
  $D$ is a divisor in a section class $e + mf$ for $m \leq l+1$, then

\begin{enumerate}  
\item $H^i (S_{k,l}, \nu ) = 0$, for $i \geq 1$.

\item $H^i (S_{k,l}, \nu \otimes \O_{S_{k,l}}(-D)) = 0$ for $i \geq 1$.
\end{enumerate}
\end{lemma}

\noindent {\bf Proof:} The line bundles $\O_{S_{k,l}}$,
$\O_{S_{k,l}}(1)$, and $\O_{S_{k,l}}(1) \otimes \O_{S_{k,l}}(-D)$ have
no higher cohomology and by Serre duality $h^2
(S_{k,l},\O_{S_{k,l}}(-D)) = 0$.  Consequently, the Euler sequence for
$\P^N$
$$0 \rightarrow \O_{\P^N} \rightarrow \oplus_1^{n+1} \O_{\P^N} (1)
\rightarrow T_{\P^N} \rightarrow 0 $$
implies that $H^i(S_{k,l},
T_{\P^N}\otimes \O_{S_{k,l}}) = 0$ and $H^i(S_{k,l}, T_{\P^N}\otimes
\O_{S_{k,l}} (-D)) = 0$ for $i \geq 1$. The standard exact sequence
\begin{displaymath}
0 \rightarrow T_{S_{k,l}} \rightarrow T_{\P^N}\otimes \O_{S_{k,l}} \rightarrow
\nu \rightarrow 0  
\end{displaymath}  
implies that $$h^i (S_{k,l}, \nu ) = h^{i+1} ( S_{k,l}, T_{S_{k,l}})
\ \ \mbox{and} \ \ h^i (S_{k,l}, \nu (-D) ) = h^{i+1} ( S_{k,l},
T_{S_{k,l}} (-D)).$$ When $i=2$ the right hand sides immediately
vanish. When $i=1$ they also vanish by Serre duality. $\Box$

\begin{theorem}\label{multiplicity1}
  When $l-k \leq 1$, the components 
\begin{equation*}
\begin{split}
  & 1. \overline{\mathcal{M}}_H (N;k,l, C, Y, I-1, J_H + 1) \\
  & 2. \overline{\mathcal{M}}_H (N-1; k,l, C,Y,I,J_H) \\
  & 3. \mathcal{X}'(\P^N; (k_0, l_0, C(k_0+l_0+1)), (k_1,l_1,
  C(k+l-k_0-l_0-1))),
\end{split}
\end{equation*}
satisfying the constraints listed in 1,2,3,5 and 6 of Theorem
\ref{main1}, occur with multiplicity one in $$D_H' \subset
\overline{\mathcal{M}} ( N;k,l,C,Y,I,J_H).$$
\end{theorem}

\noindent {\bf Proof:} By Propositions \ref{usekl} and \ref{linekl} it
suffices to restrict to the case $Y=0$, $I=1$, $\Gamma^j = H$ for all
$j$. We can then recover the general case by slicing with general
hyperplanes. 

To determine the multiplicity for the first locus we can assume that
$\Delta^I$ has dimension $N-3$. Consider the family obtained by
rotating $\Delta^I$ into $H$. Let $\Delta^I(t)$ denote this
family. We assume that $\Delta^I(0) \subset H$. Consider the family  
$$\{(t,S_{k,l}, p_1, \cdots , p_{J_H}): S \cap \Delta(t) \not= \emptyset
\}$$
of balanced scrolls $S_{k,l}$ which meet $\Delta(t)$ and have
$J_H$ marked points in $H$. The multiplicity under consideration is
the same as the multiplicity of the divisor $t=0$ in this family. The
family admits a rational map to the space of rational normal curves
with $J_H$ marked points by sending each marked surface to the
hyperplane section in $H$. This map is defined over the locus under
discussion and is smooth over that locus by Lemma \ref{normal1}.

Here we are using the fact that to show that a morphism of stacks
$\mathcal{A} \rightarrow \mathcal{B}$ where $\mathcal{B}$ is smooth
and $\mathcal{A}$ is equidimensional and smooth it suffices to show
that the Zariski tangent space to the fiber is of dimension $\dim
\mathcal{A} - \dim \mathcal{B}$.

The divisor whose multiplicity we are trying to determine occurs as a
component of the pull-back of the divisor of rational curves in
$H=\P^{N-1}$ that meet an $N-3$ dimensional linear space.  Since the
latter divisor is reduced and the morphism is smooth we conclude that
the multiplicity is one in this case.

To determine the multiplicity of the other loci it is more convenient
to look at the case when $\Delta^I = \P^N$. In that case the marking
of the point $q_I$ gives the universal surface over the loci in
question. To compute the multiplicity we can forget the marking of
$q_I$. Note that in all the loci described the hyperplane section $C$
is uniquely determined by the surface and the marked points. We get a
map from the space of surfaces to the space of rational curves by
sending the map from the surface to the map from $C$ to $H$. This is a
smooth morphism as in the previous case by Lemma \ref{normal1} and the
standard normal bundle sequence relating the normal bundle of the
reducible surface to that of the union. When $C$ is reducible, the
divisor in question is a component of the pull-back of a boundary
divisor of $\overline{M}_{0,J_H} (\P^{N-1}, k+l)$ under a smooth
morphism, hence reduced. 

Finally in case 2 and when the surface breaks into a plane union a
balanced scroll it is easy to see that the multiplicity is one by
direct computation using the determinantal representation of scrolls.
$\Box$

\begin{theorem}\label{multiplicity2}
  When $l-k \leq 1$, the components $\overline{\mathcal{M}}_{\Pi}$ and
  $\mathcal{X}_{\Pi}'$ of $$D_{\Pi}' \subset
  \overline{\mathcal{M}}_{\Pi} (N; k,l,C,Y,I,I',J_{\Pi})$$
  listed in
  $1-11$ of Theorem \ref{main2} occur with multiplicity one in
  $D_{\Pi}'$.
\end{theorem}

\noindent {\bf Proof:} It suffices to consider the case $Y=I=0$,
$I'=1$. We reduce this case to Theorem 6.2 in
\cite{vakil:rationalelliptic}. $\overline{\mathcal{M}}_{\Pi}(\P^N,
k,l,C,0,0,1,J_{\Pi})$ admits a morphism $\iota$ to
$\overline{\mathcal{M}}_{0,J_{\Pi}+1}(\P^N,d)$ which sends $(S,C,\pi)$
to $\pi: C \rightarrow \P^N$ and stabilizes. In light of \S 3, we can
interpret this morphism---at least in an open set containing the loci
we are interested in---as the morphism induced between the
corresponding Kontsevich spaces by the projection from $\F(0,1;N)$ to
$\P^N$.

We claim that $\iota$ is smooth at all the loci covered by the
theorem. Both the image and the domain of $\iota$ are equidimensional
and smooth along the loci we are interested in. It suffices to check
that the fiber is smooth. We need to compute the dimension of the
Zariski tangent space to the fiber. The Zariski tangent space to the
fiber at a point $(S,C)$, where $S$ is a Hirzebruch surface and $C$ is
a curve in a section class that lie in one of the loci covered by the
theorem, is given by $H^0 ( S, \nu_S (-C))$. If the surface is smooth,
then we can conclude that the morphism is smooth by Lemma
\ref{normal1}. If $S= S_1 \cup S_2$ has two components meeting
transversely along their common line $L$, the claim follows from Lemma
\ref{normal1} when we use the standard exact sequence for normal
bundles.

$D_{\Pi}'$ is the pull-back of the divisor $D_H$ in
\cite{vakil:rationalelliptic} by $\iota$. Since $\iota$ is smooth our
theorem follows from Theorem 6.2 in \cite{vakil:rationalelliptic}.
$\Box$ \smallskip

Finally, to conclude the multiplicity calculations we recall
Proposition 2.18 in \cite{vakil:thesis}. This proposition asserts that
if to the data of a zero dimensional locus $\mathcal{M}$, we add a
linear space of dimension $N-2$ that meets the surface or a linear
space of dimension $N-1$ that meets a marked curve, then the degree of
the stack multiplies by the degree of the surface or the degree of the
marked curve, respectively.

\section{A Simple Enumerative Consequence}

In this section we describe an application of Theorems \ref{main1} and
\ref{main2} to the characteristic numbers of balanced scrolls. We
impose enough point conditions to ensure that we can satisfy the
non-degeneracy assumptions of the theorems. When counting scrolls of
degree $n$ in $\P^{n+1}$, the hypotheses of Theorem \ref{enumerative}
require that at least $n+4$ of the linear spaces are points. It is
possible to strengthen the theorem at the expense of complicating the
algorithm.

\begin{theorem}\label{enumerative}
  Suppose $0 \leq l - k \leq 1$. Let $\{ \Delta_{a_i}^i \}_{i=1}^I$ be a
  set of linear spaces of dimension $a_i < N-2$ in $\P^N$ in general
  linear position such that

\begin{equation*}
\begin{split}
  &1. \sum_{i=1}^I N -2 - a_i \ = \ (k+l+2) \, N + 2k -4 -
  \delta_{k,l} \\
  & 2. \ \ a_I \leq N - k -l - 1 \\
  & 3. \ \  a_i = 0 \ \ \mbox{for} \ 1 \leq i \leq k+l+1 \\
  & 4. \ \ a_{k+l+2+j} \leq N-k-l-1, \ \ for \ \ \ 0 \leq j \leq a_I +
  1
\end{split}
\end{equation*}
Then there exists an algorithm which computes the number of scrolls
$S_{k,l}$ meeting $\{ \Delta_{a_i}^i \}_{i=1}^I$.
\end{theorem}

\noindent {\bf Proof:} We now describe the algorithm and prove that it
terminates without stepping outside the bounds of our non-degeneracy
assumptions. We begin with the case $N=k+l+1$ and reduce the more
general case to it later. \medskip

\noindent {\bf Step 1.} Specialize the $\Delta^i$, except for $\Delta^I$,
one by one to general linear spaces of a general hyperplane $H$ in
order of increasing dimension until a reducible solution appears.
\smallskip

In our case the first $k+l+3$ linear spaces and $\Delta^I$ are points.
We can take $H$ to be the span of the first $k+l+1$ points. We claim
that after we specialize $\Delta^{k+l+2}$ to $H$ the scrolls incident
to all $\Delta^i$ are still irreducible balanced scrolls.

The hyperplane section in $H$ has to meet the first $k+l+1$ points, so
it is non-degenerate. Similarly, since the scroll needs to meet
$\Delta^I$, it spans $\P^{k+l+1}$. By Theorem \ref{main1} if there is
a reducible scroll, then the component of the scroll in $H$ meets
$\Delta^{k+l+2}$ and contains $d+2$ of the first $k+l+1$ points, where
$d$ is the degree of the scroll. Since this is impossible the claim
follows. We repeat step 1 by specializing $\Delta^{k+l+3}$ to $H$.

Theorem \ref{main1} still applies. In this case there are 2
possibilities.

\noindent {\bf Case i.} Some scrolls can remain irreducible. Then
their hyperplane section in $H$ is the unique rational normal
curve containing the $k+l+3$ points in $H$. Repeat Step 1 by specializing
$\Delta^{k+l+4}$. Theorem \ref{main1}
still applies and this case can no longer occur. Proceed to the next
possibility. \smallskip

\noindent {\bf Case ii.} Some scrolls can become reducible. By Theorem
\ref{main1} the only reducible scrolls can be a balanced scroll of
degree $k+l-1$ in $H$ and a plane outside.

\noindent $\bullet$ If we arrived at Case ii after passing through
Case i first, this is clear since otherwise the limit hyperplane
section would be reducible. In this case the degree $k+l-1$ scroll
contains the rational normal curve and meets $\Delta^{k+l+4}$. The
plane contains $\Delta^I$. The rest of the conditions are distributed
among the two. We need to consider each way of partitioning the other
conditions in such a way that they do not impose more conditions on
either of the components than they can satisfy. (If we do not satisfy
the last clause, the algorihm will give 0.)

\noindent $\bullet$  If we are in Case ii right
after having specialized $\Delta^{k+l+3}$, then the component in $H$
needs to meet $\Delta^{k+l+3}$ and contain at least $d+2$ of the first
$k+l+2$ points if its degree is $d$. The only possibility is $d=
k+l-1$. The scroll of degree $k+l-1$ contains the $k+l+3$ points in
$H$. The plane contains $\Delta^I$. We consider each partition of the
rest of the conditions. In either case proceed to step 2.  \smallskip

\noindent {\bf Crucial Point:} Proposition \ref{direc} implies that if
a balanced scroll breaks into a union of two balanced scrolls the
gluing conditions on the directrices are automatically satisfied.
Therefore, we reduced the problem of counting scrolls of degree $k+l$
to counting pairs (Plane, Scroll of degree $k+l-1$) meeting along a
line and in addition satisfying the conditions described in Case ii.
Using Step 2, which we now describe, we further reduce the problem to
counting degree $k+l-1$ scrolls. \smallskip

Note that the scroll of degree $k+l-1$ needs to contain a line of the
plane outside $H$. The plane in turn contains a point and meets some
linear spaces.  \smallskip

\noindent {\bf Step 2.} Use Schubert calculus to reexpress the
conditions on the plane (or if one were to apply the algorithm more
generally, the conditions on the linear space containing the scroll
outside $H$) as multiples of Schubert cycles.  \smallskip

After Step 2, the plane is required to contain a point ($\Delta^I$),
meet a linear space $\Lambda_1$ in a line, and lie in a linear space
$\Lambda_2$. In turn the common line between the two scrolls is
required to meet the linear space $\Lambda_1 \cap H$ and lie in
$\Lambda_2 \cap H$. We have reduced the problem to counting degree
$k+l-1$ balanced scrolls in $\P^{k+l}$ satisfying the conditions in
Case ii) and containing a fiber lying in a linear space and meeting
another linear space. \smallskip

\noindent $\bullet$ If we arrived at this stage without going through
Case i, we are done by induction. The steps so far only used Theorem
\ref{main1} which allows for our new condition without changing the
conclusions.  In addition the scroll contains $k+l+3$ points. We can
go back to Step 1 and run the process from the beginning. \smallskip

\noindent $\bullet$ If we arrived at this stage after passing through
Case i, we have to count scrolls of degree $k+l-1$ in $\P^{k+l}$
containing a rational normal curve $C$ of degree $k+l$, meeting some
linear spaces and containing a fiber which lies in a linear space and
meets some other linear space. Proceed to step 3.  \smallskip

\noindent {\bf Step 3.}  Specialize a linear space meeting the
rational normal curve of degree $k+l$ to a general hyperplane $\Pi$ in
order of increasing dimension, but always keeping a point outside
$\Pi$, until the curve or the surface becomes reducible.

In our case Step 3 amount to breaking $C$ into a rational normal curve
of degree $k+l-1$ union a general line. Theorem \ref{main2} applies
since the surface still spans $\P^{k+l}$ and the limit of the
hyperplane section of the curve in $\Pi$ is non-degenerate. We
conclude that the surface cannot break after this degeneration. If
$\Delta^{k+l+4}$ had codimension more than 2, specialize it to $\Pi$
after which the surface has to necessarily break into a plane
containing $l$ union a degree $k+l-1$ surface. Go back and repeat
Steps 2 and 3. If $\Delta^{k+l+4}$ has codimension 2, specialize a
different linear space again in order of increasing dimension. Go back
to Step 2. We have reduced the problem to a problem of one degree
lower in $\P^{k+l-1}$.

Inductively, we reduce the problem to a problem of counting planes
with conditions of meeting a linear space, containing lines in a
linear space or containing a conic. Finally, Theorems
\ref{multiplicity1} and \ref{multiplicity2} dictate
the multiplicities with which each case occurs. This concludes the
description of the algorithm when $N= k+l+1$.  \smallskip

When $N > k + l + 1$, the algorithm is almost identical and quickly
reduces to the case $N= k+l+1$. Start by specializing $\Delta^{i}$ for
$i \leq k+l+2$ to a general hyperplane. By the assumption 3, the first
$k+l+1$ span $P$, a $\P^{k+l}$. The hyperplane section of the scrolls
in $H$ have to lie in $P$. After we specialize $\Delta^{k+l+2}$, there
are two possibilities.

{\bf 1.} If the scroll can lies in $H$, we are done by induction.

{\bf 2.} If the scroll does not lie in $H$,  $\Delta^{k+l+2}$ meets
$P$ in a point. At this stage there cannot be any reducible scrolls by
the argument given above. Specialize
$\Delta^{k+l+3}$ to $H$.

\noindent $\bullet$ If the scroll lies in $H$, we are again done by induction.

\noindent $\bullet$ If the scroll becomes reducible, proceed to Step 2
in the above process since thre must be a scroll of degree $k+l-1$ in
$P$.

\noindent $\bullet$ If the scroll remains irreducible and outside $H$,
its hyperplane section in $P$ is determined. After the next
degeneration proceed with Step 2 in the above process. This concludes
the proof.  $\Box$
\smallskip

\noindent {\bf Remark 1.} Although to satisfy the non-degeneracy assumptions
the algorithm dictates an order of specialization, the enumerative
numbers are independent of the order. By specializing the conditions
in different orders one can solve more problems. For example, to find
the number $n$ of cubic scrolls in $\P^3$ that contain a fixed twisted
cubic and five general points, we can count cubic scrolls in $\P^4$
that contain a twisted cubic, a point and meet 4 lines. If we
specialize the point to the hyperplane of the twisted cubic, some of
the limits become degenerate. Comparing the number to the one obtained
from our algorithm, we conclude that $n= 21$.  \smallskip

\noindent {\bf Remark 2.}  If we remove the non-degeneracy assumption
in Theorem \ref{main1}, there are components of $D_H$ whose general
point corresponds to a map with image a scroll $S_0$ of degree $d_0$
in $H$ with many scrolls $S_i$ outside $H$ attached to it.  The
scrolls outside $H$ can have contact of order $m_i$ with $H$ along
their common lines with $S_0$.  Moreover, the components do not have
to remain balanced. New multiplicities appear: the divisors where the
components of the scrolls have higher tangency with $H$ appear with
higher multiplicity.  The limit of the directrices usually have
tangency conditions with the limits of the hyperplane sections.

Even if we enlarge the class of problems to include these, at the next
stage worse degenerations appear. Once the surface breaks again, we
need to record the new hyperplane section which in turn can have
various tangencies with both the directrix and the old hyperplane
section. The analogue of Lemma \ref{independence} is not true for more
than two curves. When the number of curves exceeds two, I do not know
a complete list of the limits.  \smallskip

\noindent {\bf Remark 3.} We can  ask
for the characteristic numbers of $S_{k,l}$ when $l-k >1$. Theorems
\ref{main1} and \ref{main2} do not require the scrolls to be balanced.
They determine the set-theoretic limits of unbalanced scrolls. In
fact, Step 1 in the algorithm of Theorem \ref{enumerative} can be
carried out for unbalanced scrolls the same way.  However, the crucial
observation that there are no matching conditions on the directrices
of balanced scrolls no longer holds for unbalanced scrolls.  After
Step 1 of the algorithm we cannot reduce ourselves to the problem of
counting smaller degree scrolls. In addition the limit of the
hyperplane section has to meet the directrix along the special fiber.

One can reprove Theorems \ref{main1} and \ref{main2} by including this
condition. The proof is identical, only the statement and the
interpretation change. New divisors appear where the hyperplane
section contains the special fiber or the directrix thus voiding the
incidence condition. However, it becomes harder to trace this
condition during a long degeneration.

Finally, the multiplicity statements become harder for unbalanced
scrolls. Cones, especially, exhibit unexpected multiplicities.
However, in small degree one can compute the characteristic numbers of
unbalanced scrolls (see Example B2). We note that it is easy to see
that each of the degenerations in Example B2 occur with multiplicity
one by writing explicit first-order deformations, hence we omit a
detailed argument.

\section{The Gromov-Witten Invariants of $\G(1,N)$}

In this section we explain the relation between the characteristic
numbers of balanced scrolls and the Gromov-Witten invariants of
$\G(1,N)$. \smallskip

\noindent {\bf Gromov-Witten Invariants.} Recall that
$\overline{M}_{0,n}(X, \beta)$, the Kontsevich spaces of stable maps,
come equipped with n evaluation morphisms $\rho_1, \cdots, \rho_n$ to
$X$, where the $i$-th evaluation morphism takes the point $[C, p_1,
\cdots, p_n, \mu]$ to the point $\mu(p_i)$ of $X$. Given classes
$\gamma_1, \cdots, \gamma_n$ in the Chow ring $A^* X$ of $X$, the {\it
  Gromov-Witten invariant} associated to these classes is defined by 

\begin{displaymath}
I_{\beta}(\gamma_1, \cdots, \gamma_n) =
\int_{\overline{M}_{0,n}^{\text{virt}}(X,\beta)} \rho_1^*(\gamma_1)
\cup \cdots \cup 
  \rho_n^*(\gamma_n).  
\end{displaymath}

If $X$ is a homogeneous space $X= G/P$ and $\gamma_i$ are fundamental
classes of pure dimensional subvarieties $\Gamma_i$ of $X$, then there
is a close connection between the enumerative geometry of $X$ and the
Gromov-Witten invariants given by Lemma 14 in
\cite{fultonpan:kontsevich}. We reproduce this lemma for the reader's
convenience. Assume
\begin{displaymath}
\sum_{i=1}^n \mbox{codim} (\Gamma_i) = \dim(X) + \int_\beta c_1(T_X) + n-3.
\end{displaymath}
Let $g\Gamma_i$ denote the $g$ translate of $\Gamma_i$ for some $g \in
G$.

\begin{lemma}
Let $g_1 \cdots, g_n \in G$ be general elements, then the scheme
theoretic intersection
\begin{equation}\label{scheme}
\rho_1^{-1} (g_1 \Gamma_1) \cap \cdots \cap \rho_n^{-1} (g_n\Gamma_n)
\end{equation}
is a finite number of reduced points supported in $M_{0,n}(X, \beta)$
and the Gromov-Witten invariant equals the cardinality of this set
\begin{displaymath}
I_{\beta}(\gamma_1, \cdots, \gamma_n)= \# \rho_1^{-1} (g_1 \Gamma_1)
\cap \cdots \cap \rho_n^{-1} (g_n\Gamma_n). 
\end{displaymath}
\end{lemma} 
In the case of $M_{0,n}( \G(1,N),k+l)$ using Kleiman's theorem we can,
in fact, conclude that the intersection in \ref{scheme} is supported
in the locus of maps to non-degenerate curves of directrix degree
$\lfloor (k+l)/2 \rfloor$. 

Assume the $\Gamma_i$ are Schubert cycles of the form $\Sigma_{a_i}$,
the cycle of lines meeting an $a_i$ dimensional linear space.  By
Theorem \ref{scrcur} the cardinality of the intersection in
(\ref{scheme}) is equal to the number of balanced scrolls meeting
general linear spaces of dimension $a_i$, $1 \leq i \leq n$. We thus
obtain the following corollary to Theorem \ref{enumerative}:

\begin{corollary}\label{GW}
  Let $\Gamma_i = \Sigma_{a_i}$. Assume $a_i$ satisfy the conditions
  of Theorem \ref{enumerative}. Then the algorithm described in
  Theorem \ref{enumerative} provides an algorithm for computing
$$I_{k+l} (\gamma_1, \cdots, \gamma_n)$$
\end{corollary}

\noindent {\bf Remark.} One has to exercise caution when translating
the number of quadric surfaces to degree 2 Gromov-Witten invariants of
Grassmannians. Our algorithm counts actual quadric surfaces. Since 
quadric surfaces can be seen as scrolls in two distinct ways
depending on the choice of ruling, the Gromov-Witten invariant is
twice the number of quadric surfaces.  \smallskip

A closer analysis of the algorithm in Theorem \ref{enumerative} shows
that it computes the number of balanced scrolls of degree $k+l$
containing a section class of degree $k+l$ or $k+l+1$ subject to the
non-degeneracy assumptions. By an argument similar to the one just
given, our algorithm computes certain Gromov-Witten invariants of
$\F(0,1;N)$. For a sample of different approaches to the Gromov-Witten
invariants of Grassmannians and Flag manifolds see \cite{ciocan:flag},
\cite{buch:doodle} or \cite{tamvakis:doodle}.

We conclude with a table of characteristic numbers of surfaces. We use
the notation $n(N; k,l;a_0,a_1, \cdots, a_k)$ to denote the number of
$S_{k,l}$ in $\P^N$ that meet $a_0$ points, $a_1$ lines, $\cdots$,
$a_k$ $k$-planes. \bigskip

\begin{tabular}{|l|l|} \hline
\ $n(4;1,1;4,5)  =  1$ \  & \ $n(4;1,2;9,0)= 2$  \\ \hline
\ $n(4;1,1;3,7) = 9$ \  & \ $n(4;1,2;8,2) = 17$ \\ \hline
\ $n(4;1,1;2,9)= 64$ \ & \ $n(4;1,2;7,4) = 138$ \\ \hline
\ $n(4;1,1;1,11) = 430$ \ & \  $n(4;1,2;6,6) = 1140$ \\ \hline
\ $n(4;0,2;4,4) = 4$ \ & \ $n(4;1,2;5,8) = 9770$ \\ \hline
\ $n(4;0,2;3,6) = 30$ \ & \  $n(5;1,2;4,5,1) = 58$ \\ \hline
\ $n(4;0,2;2,8) = 190$ \ & \ $n(5;1,2;4,4,3) = 423$ \\ \hline
\ $n(5;1,1;3,0,8) = 48$ \ & \ $n(5;2,2;9,1) = 6$ \\ \hline
\end{tabular}

\bibliographystyle{math}
\bibliography{math}

\noindent  Mathematics Department, Harvard University, Cambridge,
  MA 02138

\noindent  E-mail: coskun@math.harvard.edu

\end{document}

%% file: lines.pstex_t
\begin{picture}(0,0)%
\includegraphics{lines.pstex}%
\end{picture}%
\setlength{\unitlength}{3947sp}%
\begingroup\makeatletter\ifx\SetFigFont\undefined%
\gdef\SetFigFont#1#2#3#4#5{%
  \reset@font\fontsize{#1}{#2pt}%
  \fontfamily{#3}\fontseries{#4}\fontshape{#5}%
  \selectfont}%
\fi\endgroup%
\begin{picture}(1449,699)(64,77)
\put(1201,164){\makebox(0,0)[lb]{\smash{\SetFigFont{12}{14.4}{\rmdefault}{\mddefault}{\updefault}{$q$}%
}}}
\put(676,239){\makebox(0,0)[lb]{\smash{\SetFigFont{12}{14.4}{\familydefault}{\mddefault}{\updefault}{$P$}%
}}}
\end{picture}